\documentclass[12pt,a4paper]{amsart}
\usepackage{amsmath, amssymb, amsthm,
amsfonts, amsthm, bm, eucal, graphicx}
\makeatletter
\@namedef{subjclassname@2020}{\textup{2020} Mathematics Subject Classification}
\makeatother
\usepackage[all]{xy}
\usepackage{graphics}
\usepackage{epsfig}
\usepackage{graphicx}
\theoremstyle{plain}
\usepackage{amssymb}
\usepackage[colorlinks]{hyperref}
\usepackage{enumerate}

\advance\hoffset-20mm \advance\textwidth40mm

\newtheorem{thm}{Theorem}[section]
\newtheorem{cor}[thm]{Corollary}
\newtheorem{lem}[thm]{Lemma}
\newtheorem{prop}[thm]{Proposition}

\theoremstyle{definition}
\newtheorem{defi}[thm]{Definition}
\newtheorem{defis}[thm]{Definitions}
\newtheorem{conj}[thm]{Conjecture}
\newtheorem{conv}[thm]{Convention}
\newtheorem{nota}[thm]{Notation}
\newtheorem{rem}[thm]{Remark}
\newtheorem{rems}[thm]{Remarks}
\newtheorem{exa}[thm]{Example}
\newtheorem{exas}[thm]{Examples}
\newtheorem{sit}[thm]{}

\newcommand{\brem}{\begin{rem}}
\newcommand{\brems}{\begin{rems}}
\newcommand{\erem}{\end{rem}}
\newcommand{\erems}{\end{rems}}
\newcommand{\bexa}{\begin{exa}}
\newcommand{\bexas}{\begin{exas}}
\newcommand{\eexa}{\end{exa}}
\newcommand{\eexas}{\end{exas}}
\newcommand{\bdefi}{\begin{defi}}
\newcommand{\edefi}{\end{defi}}
\newcommand{\bdefis}{\begin{defis}}
\newcommand{\edefis}{\end{defis}}
\newcommand{\bcor}{\begin{cor}}
\newcommand{\ecor}{\end{cor}}
\newcommand{\blem}{\begin{lem}}
\newcommand{\elem}{\end{lem}}
\newcommand{\bconv}{\begin{conv}}
\newcommand{\econv}{\end{conv}}
\newcommand{\bconj}{\begin{conj}}
\newcommand{\econj}{\end{conj}}
\newcommand{\bprop}{\begin{prop}}
\newcommand{\eprop}{\end{prop}}
\newcommand{\bthm}{\begin{thm}}
\newcommand{\ethm}{\end{thm}}
\newcommand{\bnota}{\begin{nota}}
\newcommand{\enota}{\end{nota}}
\newcommand{\bsit}{\begin{sit}}
\newcommand{\esit}{\end{sit}}
\newcommand{\be}{\begin{equation}}
\newcommand{\ee}{\end{equation}}
\newcommand{\bproof}{\begin{proof}}
\newcommand{\eproof}{\end{proof}}
\def\ba{\begin{array}}
\def\ea{\end{array}}
\def\bea{\begin{eqnarray}}
\def\eea{\end{eqnarray}}

\def\bnum{\begin{enumerate}}
\def\enum{\end{enumerate}}
\newcommand{\la}{\label}
\newtheorem*{theo*}{Theorem}

\theoremstyle{definition}

\newtheorem*{definition*}{Definition}

\def\fg{{\mathfrak g}}

\def\fl{{\mathfrak l}}

\def\fB{{\mathfrak B}}

\def\fG{{\mathfrak G}}

\def\cA{{\mathcal A}}
\def\cB{{\mathcal B}}

\def\cN{{\mathcal N}}
\def\cO{{\mathcal O}}

\def\cU{{\mathcal U}}

\def\AA{{\mathbb A}}

\def\NN{{\mathbb N}}
\def\ZZ{{\mathbb Z}}

\def\TT{{\mathbb T}}
\def\kk{{\Bbbk}}

\def\G{{\mathbb G}}
\def\Ga{{\mathbb G_{\mathrm a}}}
\def\Gm{{\mathbb G_{\mathrm m}}}

\def\PP{{\mathbb P}}

\def\ff{{\mathfrak f}}
\def\fg{{\mathfrak g}}
\def\fh{{\mathfrak h}}
\def\fii{{\mathfrak i}}
\def\fj{{\mathfrak j}}
\def\fl{{\mathfrak l}}

\def\fs{{\mathfrak s}}
\def\ft{{\mathfrak t}}
\def\fu{{\mathfrak u}}
\def\fb{{\mathfrak b}}

\def\End{\mathop{\rm End}}
\def\Aut{\mathop{\rm Aut}}
\def\Lie{\mathop{\rm Lie}}

\def\GL{\mathop{\rm GL}}

\def\SL{\mathop{\rm SL}}
\def\Jonq{\mathop{\rm Jonq}}
\def\JONQ{\mathop{\rm JONQ}}

\def\Ad{\mathop{\rm Ad}}
\def\ad{\mathop{\rm ad}}
\def\id{\mathop{\rm id}}

\def\bideg{\mathop{\rm bideg}}

\def\Norm{\mathop{\rm Norm}}
\def\Cent{\mathop{\rm Cent}}

\def\Spec{\mathop{\rm Spec}}

\def\Der{\mathop{\rm Der}}

\def\p{\partial}

\def\Aff{\mathop{\rm Aff}}

\def\Vect{\mathop{\rm Vec}}

\def\ll1{l_{\lambda}^{-1}(1)}
\def\lm1{l_{\mu}^{-1}(1)}

\def\ba{\begin{array}}
\def\ea{\end{array}}
\def\bea{\begin{eqnarray}}
\def\eea{\end{eqnarray}}

\begin{document}
\sloppy
\title[Borel subalgebras]
{Borel subalgebras of Lie
algebras of
vector fields}

\author%
{Ivan\ Arzhantsev and Mikhail\ Zaidenberg}
\address{Faculty of Computer Science, HSE University, 
Pokrovsky Boulevard 11, Moscow,
109028 Russia}
\email{arjantsev@hse.ru}
\address{Univ. Grenoble Alpes, CNRS, IF, 
38000 Grenoble, France}
\email{mikhail.zaidenberg@univ-grenoble-alpes.fr}
%
%\today
%
\begin{abstract} Let $X$ be an affine variety. 
In the present paper, we introduce 
the notion of a locally integrable Borel 
subalgebra of the Lie algebra 
 $\Lie(\Aut(X))$.
We show that these
are precisely the tangent algebras 
of the Borel subgroups of $\Aut(X)$. 
We classify the locally integrable 
Borel subalgebras in the Lie algebras 
of the groups $\Aut(S)$
of toric affine surfaces~$S$, 
notably of the affine plane
and its cyclic quotients.
See \cite{AZ25} for a description of 
the Borel subgroups 
and, more generally, 
maximal solvable subgroups
 of $\Aut(S)$.
\end{abstract}
\thanks{The work of the first author 
is supported 
by the grant RSF 25-11-00302.}
\dedicatory{To Leonid Makar-Limanov 
on occasion of his 80th birthday}
\maketitle
{\footnotesize \tableofcontents}

%%%%%%%%%%%%%%%%%%%%%%%%%%%%%%%%%%%%%%%%%%%%%%%%%%%%%%%%%%%%%%%%
\section{Introduction} 
We are working over 
an algebraically closed field $\kk$
of characteristic zero; $\AA^n$ stands 
for the affine space over $\kk$ of dimension $n$, and 
$\Ga$ resp. $\Gm$ for the additive resp. multiplicative 
group of $\kk$.

Recall that 
a \emph{Borel subalgebra}  
 of a Lie algebra $\fg$ is 
 a maximal solvable Lie subalgebra  
 of $\fg$, that is, a solvable 
Lie subalgebra that is maximal 
by inclusion among the
 solvable Lie subalgebras of $\fg$. 
 A  \emph{Borel subgroup} of an ind-group
 $\fG$ is a maximal solvable
 connected subgroup  of $\fG$. 

Let $G$ be a connected
affine algebraic group 
over $\kk$ with Lie algebra $\Lie(G)=\fg$. 
It is well known that the correspondence
$B\mapsto \Lie(B)$ yields a bijection 
between the Borel subgroups 
of $G$ and 
the Borel subalgebras of $\fg$.
 In this case, 
all Borel subgroups are conjugate, 
and all Borel subalgebras are 
$\Ad$-conjugate. The latter
is also true
for all finite dimensional 
Lie algebras  in characteristic zero:
all their Borel subalgebras are conjugate
under the action 
of the automorphism group,
see \cite[Theorem~16.4]{Hum72}.
However, this breaks down  
for some simple
Lie algebras of Cartan type 
in positive characteristic $p>0$,
as follows from 
\cite[Remark~14.3]{Hum67}. 
For example, if $p>3$, then
the first Witt Lie algebra $W_1$ 
has precisely two conjugacy 
classes of Borel subalgebras,
see \cite[Theorem~4.2]{YC15}.
See also \cite[Theorem~4.6]{BS19} 
for the conjugacy 
classes of Borel subalgebras 
of $W_n$ 
in positive characteristic $p>3$.

In the following, $X$ denotes
an irreducible affine 
algebraic variety over $\kk$. 
The automorphism 
group $\Aut(X)$ has 
a natural  affine ind-group structure, 
see e.g. \cite{FK18}.
Its Lie algebra $\Lie(\Aut(X))$ is 
a subalgebra 
of the Lie algebra of derivations 
$\Der(\cO(X))$. For example, 
the Lie algebra
$\Lie(\Aut(\AA^n))=
\Vect^{\mathrm c}(\AA^n)$ 
consists of polynomial
vector fields over $\AA^n$ 
with constant divergence, 
see e.g. \cite[Proposition~15.7.2]{FK18}, 
while 
$\Der(\kk[x_1,\ldots,x_n])=\Vect(\AA^n)$
consists of all  polynomial
vector fields over $\AA^n$.

The one to one
correspondence between 
Borel subgroups of algebraic groups
and Borel subalgebras of their 
Lie algebras
is generally not satisfied 
for $\Aut(X)$, 
even in the case $X=\AA^2$. 
Indeed, all Borel subgroups 
of $\Aut(\AA^2)$
are conjugate to the triangular subgroup 
$\JONQ(\AA^2)$, 
see \cite{BEE16}. Their Lie algebras are 
Borel subalgebras of $\Lie(\Aut(\AA^2))$;
they are pairwise conjugate 
under the adjoint action of 
$\Aut(\AA^2)$. However,  
not every Borel subalgebra
of $\Lie(\Aut(\AA^2))$ corresponds to 
a Borel subgroup of $\AA^2$, see 
Remark \ref{rem:Lie-A2}. 

It is known that the triangular subgroup 
$\JONQ(\AA^n)$ is maximal 
among the solvable 
subgroups of $\Aut(\AA^n)$
 for every $n\ge 1$, 
see \cite{FP18}.
Concerning the maximality 
of the triangular
Lie subalgebras of 
$\Lie(\Aut(\AA^n))$, we 
have the following theorem, see 
Propositions~\ref{mthm-1}
and \ref{mthm-triang}.
\bthm\la{thm:1.1}
The Lie subalgebra 
$\fj_2=\Lie(\JONQ(\AA^2))$
 of triangular 
derivations of $\kk[x,y]$
is maximal among 
the solvable subalgebras of 
$\Lie(\Aut(\AA^2))$. 
By contrast, 
the Lie subalgebra of triangular 
derivations $\fj_3=\Lie(\JONQ(\AA^3))$ 
is not maximal among 
the solvable subalgebras of 
$\Lie(\Aut(\AA^3))$. 
\ethm
Thus, 
the triangular Lie subalgebra $\fj_3$ 
is not a Borel subalgebra of 
$\Lie(\Aut(\AA^3))$. 
In Section~\ref{sec:3}
we provide a criterion for determining
whether a subalgebra 
of the Lie algebra of $\Aut(X)$ 
corresponds
to a Borel subgroup of $\Aut(X)$. 
We use the following terminology. 
\bdefis\la{def:intBsa} 
Recall that
a set $S$ of derivations of 
the structure algebra 
$\cO(X)$ is called 
\emph{locally finite} if,  
for every $f\in\cO(X)$, 
the vector subspace
$A(S)(f)\subset\cO(X)$ 
is finite dimensional, 
where $A(S)$ 
stands for the associative subalgebra 
of $\End(\cO(X))$ generated by $S$. 
Since 
$\Lie(\Aut(X))\subset \Der(\cO(X))$, 
we can speak about 
locally finite Lie subalgebras of 
$\Lie(\Aut(X))$. 

A Lie subalgebra 
$\fh$ of 
$\Der(\cO(X))$ will be called 
\emph{integrable}
if $\fh\subset\Lie(G)$
for a connected algebraic subgroup 
$G\subset\Aut(X)$. 
By \cite[Theorem E]{KZ24},
$\fh$ is integrable if 
and only if it is locally finite; 
in the latter case,
$\fh\subset \Lie(\Aut(X))$.

A Lie subalgebra $\fh$  of $\Lie(\Aut(X))$ 
will be called 
\emph{locally integrable}
if it is filtered by a sequence of 
locally finite Lie subalgebras.
If a solvable $\fh$ is locally integrable 
and maximal among the solvable 
 locally integrable subalgebras of $\Lie(\Aut(X))$, 
then we say that $\fh$ is a 
\emph{locally integrable Borel subalgebra} 
 of $\Lie(\Aut(X))$. 
\edefis
Note that a Lie subalgebra
$\fh\subset\Lie(\Aut(X))$
 is locally integrable
if and only if $\fh\subset \Lie(G)$, 
where $G$
is a connected nested subgroup 
of $\Aut(X)$, 
cf. \cite{KZ24}.
The following result is central 
to our paper,
see Corollary \ref{cor:corresp}.
\bthm\la{thm:1.2} The correspondence 
$G \mapsto \Lie(G)$
leads to a one to one correspondence 
between the conjugacy classes of 
the Borel subgroups 
of $\Aut(X)$ 
and
the $\Ad$-conjugacy classes 
of the locally 
integrable 
Borel subalgebras
of $\Lie(\Aut(X))$.
\ethm
The proof explores several recent results
concerning the structure of
unipotent and solvable subgroups 
of $\Aut(X)$, see \cite{CKRVS26,
Kra25, KZ24, Per23}; 
see also Theorem \ref{thm:nested} and  
Lemma \ref{lem:nested} below. 

Recall that
every Borel subgroup of $\Aut(\AA^2)$
 is nested and conjugate to 
the triangular subgroup 
$\JONQ(\AA^2)$ (\cite{BEE16}). 
The same applies to Borel subgroups 
of the group $\mathfrak{U}_2$ 
 of automorphisms of $\AA^2$ 
 whose Jacobian 
is equal to 1 (\cite{BEE16}) and 
 to Borel subgroups 
of the group 
$\Aut_0(\AA^2)$ of automorphisms
 fixing the origin
(\cite{AZ25}). 

As follows from Theorem \ref{thm:1.2},
any locally integrable Borel subalgebra of
 $\Lie(\Aut(\AA^2))$ is 
 $\Ad$-conjugate to the 
triangular  Borel subalgebra $\fj_2$ 
of $\Lie(\Aut(\AA^2))$ (which is maximal
among the solvable subalgebras of 
$\Lie(\Aut(\AA^2))$,
see Theorem \ref{thm:1.1}). 

The group $\Aut(\AA^3)$ has non-conjugate 
Borel subgroups, see \cite{RUVS25}.
Hence, by Theorem~\ref{thm:1.2},
it also has locally integrable Borel subalgebras 
that are not Ad-conjugate to
the subalgebra 
${\fj_3\subset\Lie(\Aut(\AA^3))}$
of the triangular derivations. 

Given a primitive root of unity $\zeta$ 
of order $d>1$ and 
an integer number $e$ with
 $1\le e<d$ and $\gcd(d,e)=1$, 
 we denote by
$G_{d,e}\subset\Aut(\AA^2)$
the cyclic group  generated by 
$(x,y)\mapsto (\zeta^e x, \zeta y)$.
The quotient 
$X_{d,e}=\AA^2/G_{d,e}$ 
is a normal affine toric surface.
Any normal affine toric surface
is isomorphic to some $X_{d,e}$
except for $\AA^2$, 
$\AA^1\times\AA^1_*$ and 
$(\AA^1_*)^2$,
where 
$\AA^1_*:= \AA^1\setminus\{0\}$. 

The groups $\Aut(X_{d,e})$ 
are amalgams, see 
\cite{AZ13}. Any Borel subgroup 
 of $\Aut(X_{d,e})$ is isomorphic 
 to a triangular one.
 However, the lower and 
upper triangular subgroups
 of $\Aut(X_{d,e})$
 are conjugate if and only if 
 $e^2\equiv 1\mod d$. 
Otherwise, there exist exactly
two conjugacy classes of 
Borel subgroups 
of $\Aut(X_{d,e})$,  
see \cite[Theorem~1.1]{AZ25} 
or Theorem \ref{main} below. 
The next corollary follows
immediately
from this last result 
and from Theorem \ref{thm:1.2}.
\bcor
If $e^2\equiv 1\mod d$, then
there exists only one $\Ad$-conjugacy class 
of locally integrable Borel subalgebras of 
$\Lie(\Aut(X_{d,e}))$.
Otherwise, there exist exactly
two such classes. 
\ecor
Recall that the 
Lie algebra $\Der(X)$ 
is simple if and only if $X$ is smooth, 
see \cite{Jor86} and \cite{Sie92};
see also  \cite{BF18} and
\cite{HM94} 
for alternative proofs.
The Lie algebra 
$\Vect^0(\AA^n)$ of zero divergence 
vector fields on $\AA^n$ also is simple and is 
the unique proper ideal of the Lie algebra
$\Vect^{\mathrm c}(\AA^n)=\Lie(\Aut(\AA^n))$  
of constant divergence vector fields, 
see \cite[Lemma~3]{Sha81}.
For the affine toric surfaces $X_{d,e}$ 
(which are singular for $d>1$), we 
establish the following results.

\bthm\la{thm:1.3} $\,$

\bnum\item[$(a)$]
The Lie algebra 
$\Lie(\Aut(X_{d,e}))$ 
can be naturally identified with 
a subalgebra of 
$\Lie(\Aut(\AA^2))^{\Ad(G_{d,e})}$.
\item[$(b)$] 
Suppose that $e|(d+1)$. Then 
$\Lie(\Aut(X_{d,e}))=
\Lie(\Aut(\AA^2))^{\Ad(G_{d,e})}$.
\item[$(c)$]  
The subalgebra $\fg_{d,e}\subset
\Lie(\Aut(X_{d,e}))$
of zero divergence vector fields 
is not simple.
\enum
\ethm
See Lemma
 \ref{lem:inclu} and 
Corollary \ref{cor:5.8} for (a), 
 Example \ref{exa:veron} 
 (where $e=1$)
and Proposition \ref{lem:5.13} 
(where $e>1$)
for (b), and Proposition 
\ref{prop:nonsimple}  for (c). 
We do not know whether 
(b) holds for all pairs $(d,e)$. 

The content of this paper is as follows.
Section \ref{sec:2} contains the necessary 
preliminaries.
Theorem \ref{thm:1.2} is proved
 in Section \ref{sec:3}, 
Theorem \ref{thm:1.1} in 
 Section \ref{sec:4}, and Theorem \ref{thm:1.3}
  in Section \ref{sec:5}. 
\section{Preliminaries}\la{sec:2}
\subsection{Tangent Lie algebras 
of derivations}
\bdefi[{\rm cf. \cite[Definition~4.1.1]{KZ24}}]
A Lie algebra is
called
\emph{nested} if it has an ascending filtration 
by a sequence of finite-dimensional 
Lie subalgebras. 
\edefi
For an affine variety $X$, 
the automorphism group $\Aut(X)$ 
is an affine ind-group, 
see e.g.~\cite{FK18}. 
The Lie algebra of derivations $\Der(\cO(X))$ 
coincides with 
the Lie algebra of vector fields $\Vect(X)$.
The tangent Lie algebra 
$\Lie(\Aut(X))$
is a subalgebra of $\Der(\cO(X))$. 
For example, $\Lie(\Aut(\AA^n))$
consists of derivations 
 in $\Der(\kk[x_1,\ldots,x_n])$ 
 with constant divergence, 
 see \cite[Proposition~15.7.2]{FK18}.
For $n\ge 2$, $\Lie(\Aut(\AA^n))$ 
is a maximal proper Lie subalgebra of 
 $\Der(\kk[x_1,\ldots,x_n])$, 
 see \cite[Proposition 2.21]{Bav17}.

Note that the Lie algebra $\Lie(\Aut(X))$ 
viewed as a vector space
is of countable dimension. 
Indeed, let us choose
a dense open subset $U\subset X$ 
such that $TX|_U$ 
is a trivial vector bundle of rank 
$n=\dim(X)$. 
The vector space $\Vect(X)$
is clearly a subspace of 
$\Vect(U)\simeq \cO(U)^n$. 
Since $\cO(U)$ is of countable dimension, 
the same is true for $\Vect(X)$ and
$\Lie(\Aut(X))$.

The following theorem gives 
a uniform estimate 
for the derived lengths 
of solvable Lie subalgebras 
of derivations, 
cf. also  \cite[Theorems~1.1--1.3]{MR14}. 
\bthm[{\rm \cite[Corollary~4]{MP14}}]\la{thm:MP}
Let $L$ be a nilpotent (solvable, respectively) 
subalgebra of $\Der(\cO(X))$. 
Then the derived length of $L$ is at most $n$ 
($2n$, respectively).
\ethm
We need the following simple lemma.
Its second part follows from Proposition~3.10 in 
\cite{FP18} and its proof.
\blem\la{lem:solv-union} Let $\fh$ be a Lie algebra. 
An ascending chain of solvable subalgebras 
$\fh_\alpha$ of  $\fh$,
whose derived lengths are uniformly bounded,
has a solvable union whose derived length 
satisfies the same bound. 
An analogous statement holds for 
an ascending chain of solvable subgroups 
of a given group, whose derived lengths
are uniformly bounded. 
\elem
\bproof We present the proof for the case of
a Lie algebra; it adapts mutatis mutandis 
to the case of a group. We write $\fh^{(0)}=\fh$,
$\fh^{(r+1)}=[\fh^{(r)},\,\fh^{(r)}]$.

By our assumption, $\fh_\alpha^{(d)}=0$ 
for some $d>0$ 
and for every index $\alpha$. 
Let $\fh_0=\bigcup_\alpha \fh_\alpha$. 
Any element $w\in \fh_0^{(d)}$
is a Lie word formed from a finite set of elements 
$h_1,\ldots,h_n$ in $\fh_0$.
Choose $\alpha$ such that $h_i\in \fh_\alpha$ 
for $i=1,\ldots,n$.
Then $w\in\fh_\alpha^{(d)}=0$, which implies 
that $\fh_0^{(d)}=0$. 
\eproof
\bcor\la{cor:maximal}
Any solvable Lie subalgebra 
$L\subset\Der(\cO(X))$ 
is contained in a Borel subalgebra.
The same is true for 
the Lie subalgebras of $\Lie(\Aut(X))$.
\ecor
\bproof This follows 
from Zorn's Lemma 
due to Lemma \ref{lem:solv-union} 
and the uniform boundedness 
of derived lengths in 
totally ordered chains of Lie subalgebras,
see Theorem \ref{thm:MP}.
\eproof
\brem\la{rem:sl2} The Lie algebra $\Der(\cO(X))$ 
may be unsolvable even 
if $\Aut(X)$ is a solvable algebraic group. 
For example, this  is the case for 
$X=\AA^1$. Indeed,
$\Aff(\AA^1)=\Aut(\AA^1)\cong\Gm\ltimes\Ga$
 is a metabelian affine algebraic
group, while $\Der(\cO(\AA^1))$ contains 
a simple Lie subalgebra $\fh\simeq {\rm \fs\fl}(2,\kk)$ 
generated by 
\[e=x^2 \p/\p x,\quad f=-\p/\p x,\quad\text{and}\quad
h=2x\p/\p x,\]
see \cite[Proposition~3]{AMP11}. 
The solvable Lie subalgebra
\[\Lie(\Aut(\AA^1))=\kk\p/\p x\oplus\kk x\p/\p x
\subset\Der(\cO(\AA^1))\] 
consists of vector fields on $\AA^1$ 
with constant divergence. 
\erem
\subsection{Correspondence 
``ind-subgroups 
$\leftrightsquigarrow$ 
tangent subalgebras''}
In this subsection we mainly 
follow \cite{KZ24}. 
We begin with the following remarks. 
\brems\la{rem:loc-fin}  $\,$

1. A Lie subalgebra 
$\fh\subset \Der(\cO(X))$ 
is locally finite (see Definitions \ref{def:intBsa}
in the Introduction)
if and only if every finite 
dimensional subspace 
$V\subset \cO(X)$ is contained 
in an $\fh$-invariant  
finite dimensional subspace 
$\tilde V\subset \cO(X)$.
A locally finite Lie subalgebra 
$\fh\subset \Der(\cO(X))$ 
is finite dimensional, see 
\cite[Lemma~1.6.2]{KZ24}. Indeed, 
let us choose an $\fh$-invariant 
finite dimensional 
subspace  $E\subset\cO(X)$ 
that contains a system of
 generators of the algebra
$\cO(X)$. Then the  linear 
restriction map
$\fh\to \fh|_E\subset\End(E)$
 is injective. 

2.
Let $G\subset \Aut(X)$ be 
an algebraic subgroup.
Then $\Lie(G) \subset\Lie(\Aut(X))$ 
is a locally finite Lie subalgebra. 
However, a locally finite
 Lie subalgebra of $\Lie(\Aut(X))$ 
 is not necessarily
 the tangent algebra of an 
 algebraic subgroup 
 of $\Aut(X)$.
 
If $G\subset \Aut(X)$ is a nested 
ind-subgroup,
then $\Lie(G) \subset\Lie(\Aut(X))$ is 
a nested Lie algebra
filtered by locally finite Lie subalgebras.
A partial converse of this is also true.
\erems
\bprop[{\rm \cite[Proposition~4.2.1]
{KZ24}\footnote{The adjective "unipotent" 
is misleading in part (1) of
\cite[Proposition~4.2.1]{KZ24}.
By deleting it, the proposition
and its proof remain valid.}}] 
\la{prop:KZ}
Let $\fh= \varinjlim \fh_i \subset 
\Lie(\Aut(X))$
be a nested Lie subalgebra.
Then the following hold.
\bnum
\item[$(1)$] If $\fh$ is generated by 
locally nilpotent derivations,
then there exists a nested subgroup
$G= \varinjlim G_i \subset\Aut(X)$ 
such that 
$\fh= \Lie(G)=\varinjlim\Lie(G_i)$.

\item[$(2)$] Suppose that all $\fh_i $ 
are locally finite.
Then there exists a unique 
smallest connected nested ind-subgroup 
$G=G_{\min}(\fh) \subset \Aut(X)$
such that $\fh \subset \fg:=\Lie(G)$.
Moreover, a closed ind-subgroup 
$\hat G \subset \Aut(X)$ contains $G$
if and only if $\fh \subset \Lie(\hat G)$.
\enum
\eprop
The following is a version of
\cite[Theorem~5.1.1]{KZ24}. \footnote{
In the original formulation, Theorem 5.1.1
 in \cite{KZ24} is
applied to an arbitrary affine ind-group, 
not necessarily 
of the form $\Aut(X)$.
The proof of this theorem in \cite{KZ24} 
is based on Lemma~5.1.6 in \cite{KZ24}. 
However, this lemma turns our to be false, 
and so Theorem 5.1.1
therefore requires a proof.  In 
the present restricted version, 
Theorem \ref{thm:KZ0}
is due to Hanspeter Kraft and 
the second author. }
 \bthm \la{thm:KZ0}
Let $G\subset\Aut(X)$ 
be a solvable subgroup generated 
by a family of
connected algebraic subgroups 
$H_i\subset G$, $i\in I$.
\bnum
\item[$(1)$]
If $I$ is finite, then $G$ is 
a solvable algebraic subgroup 
with Lie algebra
$\Lie(G) = \langle\Lie(H_i)\, |\,  
i \in I\rangle_{\Lie}$.
\item[$(2)$]
In general, $G$ is nested and is of the form 
$G = U_G \ltimes T$, where 
$T \subset \Aut(X)$ is a torus
and $U_G$ is a nested unipotent group.
\item[$(3)$] If $G$ is generated 
by a family of 
unipotent algebraic subgroups, then 
$G = U_G$ is a nested unipotent group.
\enum
\ethm
The proof relies on the following recent 
and important result, kindly communicated 
by Hanspeter Kraft to the second author.
It confirms a conjecture in
\cite[Conjecture~6.6]{KZ25}.
See also an earlier version 
\cite[Theorem~B]{CRX23} 
concerning 
connected abelian subgroups 
of $\Aut(X)$.
\bthm[{\rm \cite[Theorem~A]{CKRVS26};
see also \cite[Theorem~5]{Kra25}}] \la{thm1}
Let $X$ be a quasi-affine algebraic variety,
$A$ be an irreducible algebraic variety, 
and 
$f\colon A\to \Aut(X)$ be a morphism 
such that ${\rm id}_X \in f(A)$.
Let $G$ be the subgroup of $\Aut(X)$ 
generated by the image $f(A)$.
 If $G$ is solvable,
then $G$ is a connected algebraic 
subgroup of $\Aut(X)$. 
\ethm
\bproof[Proof of Theorem \ref{thm:KZ0}] 
For the implications (1)$\Longrightarrow$(2) 
and 
(2)$\Longrightarrow$(3), see 
the proof of Theorem 5.1.1 in \cite{KZ24}. 
To prove (1),
we note that 
$G=\langle H_1,\ldots,H_n\rangle$ 
is generated by the product $H_1\cdots H_n$. 
Let, in Theorem \ref{thm1},
$A=H_1\times\cdots\times H_n$ and 
$f\colon A\to \Aut(X)$ be 
the multiplication map.
By this theorem, we conclude that 
$G=\langle f(A)\rangle\subset\Aut(X)$
is an algebraic subgroup.
\eproof 
\bcor\la{thm:KZ2} 
Any solvable subgroup $G\subset\Aut(X)$ 
generated by 
connected algebraic subgroups is 
nested of derived length 
$\le\dim(X)+1$. 
\ecor
\bproof 
By Theorem  \ref{thm:KZ0}(2),
$G$ is nested of the form 
$G = U_G \ltimes T$.
A nested unipotent subgroup 
$U\subset \Aut(X)$ 
is solvable with derived length 
$\le\max\{\dim(Ux) \,| 
\,x\in X\} \le \dim(X)$, 
see \cite[Theorem~5.3.1]{KZ24}.
Therefore,  $G$
is of derived length
$\le\dim(X)+1$.
\eproof
\bdefi
Abusing the language, we say that 
an ind-group $G$ is \emph{nested 
by connected algebraic subgroups},
if it admits a countable ascending 
filtration 
by connected algebraic subgroups.  
\edefi
\brems $\,$

\begin{itemize}\item
 If a subgroup $G\subset\Aut(X)$ is nested 
by connected algebraic subgroups, then 
it is algebraically generated in the sense 
of \cite{KZ24} and $\Lie(G)$ is locally integrable.
\item
If $\fg\subset\Lie(\Aut(X))$ is 
a locally integrable subalgebra, 
then it consists of locally finite elements.
\end{itemize}
\erems

Corollary \ref{thm:KZ2} can be strengthen 
as follows.
\bthm[{\rm\cite[Theorem~A]{CKRVS26};  
see also~\cite[P.~29, Corollary~B]{Kra25}}]
\la{thm:nested} Let $G\subset\Aut(X)$ 
be a
solvable closed connected subgroup. 
Then $G$ is nested by a sequence of
connected algebraic subgroups $G_i$
and has a derived length  $\le\dim(X)+1$. 
The Lie algebra $\fg=\Lie(G)$
is locally integrable and 
solvable with derived length $\le\dim(X)+1$.
\ethm
For the reader's convenience, 
we provide an argument. 
\bproof The closed ind-subgroup 
$G\subset\Aut(X)$ 
has an ascending filtration 
by irreducible algebraic subsets $A_i$
equivalent to the original filtration, 
see
\cite[Proposition~1.6.3]{FK18}. 
We can replace any 
$A_i$ by $A'_i$ such that $A_i'$ is 
an  irreducible algebraic subset containing
$\id_X$ and the $A'_i$ 
form an equivalent filtration again. 
Indeed, let us take 
an irreducible curve $C\subset G$
that contains $\id_X$ and an element 
$a_i^{-1}$
where $a_i\in A_i$. 
The image $C A_i$ of $C\times A_i$
under the multiplication morphism 
$C\times A_i\to \Aut(X)$, $(c,a)\mapsto ca$, 
contains $A_i$
and is an irreducible constructible subset 
of some $A_{n(i)}$. 
Now replace $A_i$ by
 the closure $A'_i=\overline{CA_i}$. Then
$G_i=\langle A'_i\rangle$
is a solvable algebraic subgroup of $G$, 
see Theorem \ref{thm1}.
Therefore, $G=\varinjlim G_i$ is nested 
by the connected 
algebraic subgroups $G_i$.  
By Corollary \ref{thm:KZ2}, 
the derived length of $G$ is $\le\dim(X)+1$.
This proves the first assertion.
Now, the second assertion follows by applying
 \cite[Lemma~5.1.4(2)]{KZ24} and its proof. 
\eproof
Let us also mention the following theorem.
\bthm[{\rm\cite[Theorem~1.4]{Per23}}]
\la{thm:closed} 
Any connected nested subgroup 
$G \subset\Aut(X)$ is closed.
\ethm
\brem\la{rem:ques}
A connected ind-group $\fG$ 
is commutative if $\Lie(\fG)$ 
is, see \cite[Corollary~7.5.3]{FK18}. However,
 we do not know whether a 
 connected ind-group $\fG$ 
is solvable if $\Lie(\fG)$ is, 
see \cite[Question~5]{KZ24}. 
\erem
The following partial result is implicit 
in the proof of 
\cite[Proposition~4.2.1]{KZ24} 
(see Proposition~\ref{prop:KZ}). 
For the reader's convenience, 
we provide a proof. 
\bprop\la{prop:solv}
Let 
$\fh= \varinjlim \fh_i \subset \Lie(\Aut(X))$
be a solvable locally integrable Lie subalgebra, 
where all $\fh_i$ 
are locally finite 
(see Definitions \ref{def:intBsa}). 
Let $G=G_{\min}(\fh) \subset \Aut(X)$ 
be 
the unique smallest connected 
nested ind-subgroup
such that $\fh \subset \Lie(G) $, 
see Proposition \ref{prop:KZ}(2). 
Then the following hold.
\bnum\item[$(1)$] 
$G$ is solvable and closed 
with derived length 
$\le\dim(X)+1$. 
\item[$(2)$] If $\fh$ is generated 
by locally nilpotent derivations, 
then it consists of locally nilpotent 
derivations, the group $G$ 
is unipotent and 
$\fh=\Lie(G)$. 
\enum
\eprop
Before moving on to the proof, 
let us first 
recall 
some definitions and
Lemmas \ref{lem:1.8.2}--\ref{lem:3.2.1} 
borrowed in \cite{KZ24},
 then deduce an additional 
 Lemma \ref{lem:new}. 
\bdefis $\,$

\begin{itemize} \item
A finite-dimensional Lie subalgebra 
 $\ft\subset\Der(\cO(X))$  
 is called \emph{toral} 
 if it consists 
 of semisimple  (locally finite) 
 elements. 
 
\item A locally finite Lie subalgebra 
$\fg \subset \Der(\cO(X))$ 
is called 
\emph{$J$-saturated} if 
 it contains, with each element $\mu$, 
 the semisimple 
 and the locally nilpotent parts 
 $\mu_{\mathrm s}$ and 
 $\mu_{\mathrm n}$.
 \end{itemize}
 \edefis
 A toral subalgebra $\ft$ 
is abelian 
 and therefore locally finite, 
 see \cite[Section~8.1]{Hum72} 
and \cite[Section~1.8]{KZ24}. 
 \blem[{\rm\cite[Lemma~1.8.2]{KZ24}}]\la{lem:1.8.2}
For a toral Lie subalgebra 
$\ft\subset \Der(\cO(X))$, 
there exists a unique  smallest torus 
$T\subset\Aut(X)$ such that 
$\ft\subset\Lie(T)$. 
Furthermore, if a subspace 
$E\subset \Der(\cO(X))$ 
is invariant under the adjoint 
action of $\ft$, 
then $E$ is invariant under $T$.
 \elem
\blem[{\rm see 
\cite[Lemmas~1.6.3 and~3.1.3]
{KZ24}}]\la{lem:3.1.3}
Let $\fg \subset \Der(\cO(X))$ be 
a solvable locally finite  
Lie subalgebra and $\ft$ 
be a maximal toral subalgebra of $\fg$. 
Suppose $\fg$ is $J$-saturated. 
Then the locally nilpotent elements 
from $\fg$ form 
a Lie ideal $\fg_{\mathrm n}$ of $\fg$, 
and 
$ \fg=\ft \ltimes\fg_{\mathrm n}$. 
\elem
\blem[{\rm\cite[Lemma~3.2.1(1)]
{KZ24}}]\la{lem:3.2.1}
Let $\fh\subset \Vect(X)$ be 
a solvable locally finite 
 Lie subalgebra. 
Then, the smallest $J$-saturated 
Lie subalgebra 
$\fg=\fg_{\min}(\fh)\subset \Vect(X)$ 
containing $\fh$  
is solvable and locally finite.
\elem
In the proof of Proposition \ref{prop:solv}, 
we use the following lemma 
(cf. the proof of
 \cite[Proposition~4.2.1]{KZ24}).
\blem\la{lem:new}
Let $\fh\subset\Lie(\Aut(X))$ be 
a  locally finite solvable 
Lie subalgebra. 
Then there exists a unique 
minimal solvable connected 
algebraic subgroup 
$G=G_{\min}(\fh)\subset \Aut(X)$ 
with derived length 
$\le \dim(X)+1$ such that 
$\fh\subset\Lie(G)$. 
\elem
\bproof
Let $\fg=\fg_{\min}(\fh)$ be
the smallest locally finite
solvable
$J$-saturated Lie subalgebra 
of $\Der(\cO(X))$ that contains $\fh$, 
see Lemma \ref{lem:3.2.1}.  
We have a semidirect product 
decomposition
$\fg= \ft
\ltimes\fg_{\mathrm n}$,
where $\ft$
 is a maximal toral Lie 
 subalgebra of $\fg$
and $\fg_{\mathrm n}$ is 
a nilpotent ideal
consisting of locally 
nilpotent elements, 
see Lemma \ref{lem:3.1.3}. 
Moreover, 
there exists a minimal torus 
$T \subset\Aut(X)$ 
with $\ft\subset\Lie(T)$
such that $T$ normalizes 
$\fg_{\mathrm n}$, 
see Lemma \ref{lem:1.8.2}. 
The algebraic subgroup 
$U:=\exp(\fg_{\mathrm n})\subset\Aut(X)$
is unipotent with Lie algebra 
$\Lie(U) =\fg_{\mathrm n}$. 
The product $G= TU=
U T \simeq T\ltimes U$
is a solvable connected algebraic 
subgroup of $ \Aut(X)$
such that $\fh\subset\fg\subset\Lie(G)$.
By construction, $G$ is 
the smallest subgroup 
possessing this property. 
By Corollary \ref{thm:KZ2},
the derived length of $G$ 
is bounded above by
$\dim(X)+1$. 
\eproof
\bproof[Proof of  
Proposition \ref{prop:solv}]
By Lemmas \ref{lem:3.2.1} 
and \ref{lem:new}, 
we have 
$\fh_i\subset \fg_{\min}(\fh_i)\subset 
\Lie(G_{\min}(\fh_i))$.
By the uniqueness of  subalgebras 
$\fg_i:= \fg_{\min}(\fh_i)$ 
and solvable connected algebraic subgroups 
$G_i:=G_{\min}(\fh_i)$ 
we conclude that 
$\fg_i\subset\fg_{i+1}$ and 
$G_i\subset G_{i+1}$.  

Let $\hat G=\varinjlim(G_i)$.
Certainly, $\hat G$ is 
a solvable connected 
nested subgroup of $\Aut(X)$
with derived length  $\le \dim(X)+1$ 
and
$\fh\subset \Lie(\hat G)=\varinjlim \Lie(G_i)$.
Since $G$ and $\hat G$ 
are connected 
and nested, 
they are closed, 
see Theorem \ref{thm:closed}.
By Proposition \ref{prop:KZ}(2), we have 
$G=G_{\min}(\fh)\subset\hat G$.
Therefore, $G$ is solvable and closed
with derived length  $\le \dim(X)+1$.
This proves part (1).

To prove part (2), suppose that 
$\fh$ is generated by locally nilpotent 
derivations $\p_\alpha$, $\alpha\in A$. 
By Proposition \ref{prop:KZ}(1), there 
exists a nested subgroup
$G'= \varinjlim G'_i \subset\Aut(X)$ 
such that 
$\fh= \varinjlim\Lie(G'_i)$. 
We can assume that 
$G'$ and all $G'_i$ are connected. 
Since all $\Lie(G'_i)$ are solvable, 
the connected 
algebraic subgroups 
$G'_i$ are also solvable. 

Now consider the unipotent 
one-parameter 
subgroups 
$U_\alpha=\exp(\kk\p_\alpha)$. 
Since  $\p_\alpha\in\Lie(G'_i)$
for some $i=i(\alpha)$, we have 
$U_\alpha\subset G'_i$. For a finite set 
of indices 
$\eta=\{\alpha_1,\ldots,\alpha_n\}\subset A$ 
we can find 
$i=i(\eta)$ such that 
$\cU_\eta:=\langle U_{\alpha_1},
\ldots,U_{\alpha_n}\rangle\subset G'_i$.
It follows that $\cU_\eta$ is a solvable 
connected affine algebraic group 
generated by 
unipotent subgroups. Certainly, 
such a group is unipotent. 

On the other hand, since by our assumption 
$\fh=\langle\p_\alpha\,|\,\alpha\in A\rangle_{\Lie}$,
for every 
$i$ we can find 
$\eta_i=\{\alpha^{(i)}_1,
\ldots,\alpha^{(i)}_{n_i}\}\subset A$
such that 
\[\Lie(G'_i)\subset\langle
\p_{\alpha^{(i)}_1},
\ldots,\p_{\alpha^{(i)}_{n_i}}\rangle_{\Lie}.\]
It follows that $G_i'\subset
U_{\eta_i}$ (see 
\cite[Remark~17.3.3]{FK18}) is 
a unipotent algebraic group. 
Thus, $G'=\varinjlim G'_i
=\varinjlim U_{\eta_i}$
is a nested unipotent ind-group, and 
$\fh=\varinjlim \Lie(U_{\eta_i})$ consists
of locally nilpotent derivations.

By Proposition \ref{prop:KZ}(2),
we have $G\subset G'$.  Hence, 
$G$ is also a nested unipotent ind-group. 
As shown, 
$\fh\subset \Lie(G)\subset\Lie(G')=\fh$,
therefore $\Lie(G)=\fh$, as stated. 
\eproof
The following question arises naturally 
in view of Theorem \ref{thm:KZ0}(1)
and Propositions \ref{prop:KZ} 
and~\ref{prop:solv}. Cf. also closely related 
Questions 2 and 3 in \cite{KZ24}.

\medskip

\noindent {\bf Question 1.}\emph{
Let 
$\fh=\langle \fh_1,\ldots,\fh_n\rangle_{\Lie}$
be a solvable Lie subalgebra of 
$\Lie(\Aut(X))$ generated by 
the locally finite Lie subalgebras 
$\fh_i\subset \Lie(\Aut(X))$. 
Is it true that $\fh$ is locally finite, 
and therefore also finite dimensional?} 

\medskip

See \cite{Zai26} for some partial results.

The following example shows that 
the local finiteness assumption 
in Question 1
is essential.
\bexa Let 
$X=\AA^1\times\AA^1_*=
\Spec\kk[x,y,y^{-1}]$.
Consider the metabelian Lie subalgebra 
$\fh\subset
\Der(\kk[x,y,y^{-1}])$ generated by
 the derivations
\[\delta=\p/\p y\quad\text{and}\quad \p = 
y^{-1} \p /\p x.\] 
The vector space $\fh$ is infinite 
dimensional with a basis 
$\delta,\p_1=\p,\p_2,\ldots$, 
where $\p_i=y^{-i} \p /\p x$.
Notice that $\delta\notin\Lie(\Aut(X))$,
and it is not locally finite. 
\eexa
\section{Correspondence 
``Borel subgroups 
$\leftrightsquigarrow$ 
Borel subalgebras''}\la{sec:3}
In this section we prove Theorem \ref{thm:1.2}.
Let us start with the following simple lemma.
\blem\la{lem:nested}
Every Borel subgroup $\cB\subset\Aut(X)$ is a
closed subgroup of derived length  
${\le\dim(X)+1}$
nested by connected algebraic subgroups. 
\elem
\bproof The closure of a solvable subgroup 
 is also a solvable subgroup, see 
\cite[Lemma 5.1.3(2)]{KZ24}. Therefore,
$\cB$ is closed.  
Theorem \ref{thm:nested} now
applies and gives the result. 
\eproof
We have the following analogue of Corollary 
\ref{cor:maximal} 
 (see Definition \ref{def:intBsa}).
\blem Any solvable locally 
integrable subalgebra 
$\fh$
of $\Lie(\Aut(X))$ 
is contained in a locally integrable 
Borel subalgebra
$\fb\subset\Lie(\Aut(X))$.
\elem
\bproof
Let $\{\fh_\alpha\,|\,\alpha\in A\}$ be a 
totally ordered ascending chain of 
Lie subalgebras of $\Lie(\Aut(X))$
that contain $\fh$ and are filtered 
by solvable locally finite subalgebras.
Let 
$G_\alpha=G_{\min}(\fh_\alpha)$, 
see Proposition \ref{prop:solv}.
Then $\{G_\alpha\,|\,\alpha\in A\}$ is a
totally ordered ascending chain of solvable 
connected nested subgroups of $\Aut(X)$. 
By Lemma \ref{lem:solv-union} and
Theorem \ref{thm:nested}, 
the closure $\bar G$ of 
$G:=\varinjlim G_\alpha$
is a solvable 
connected subgroup of $\Aut(X)$
of derived length $\le \dim(X)+1$
filtered by an ascending sequence 
of connected algebraic subgroups
$\bar G_i$. Then 
$\fg:=\Lie(\bar G)=\varinjlim \Lie(\bar G_i)
\subset\Lie(\Aut(X))$ is 
a solvable Lie subalgebra filtered by 
locally finite Lie subalgebras $ \Lie(\bar G_i)$.
Since $\fh_\alpha\subset\Lie(\bar G)$
for all $\alpha\in A$, $\Lie(\bar G)$
is an upper bound for the chain 
$\{\fh_\alpha\}$.
By Zorn's Lemma, the set of 
solvable Lie subalgebras of $\Lie(\Aut(X))$
that contain $\fh$ and are filtered 
by locally finite subalgebras has 
a maximal element $\fb$. 
Thus, $\fb$ is a locally integrable 
Borel subalgebra which contains $\fh$.
\eproof
\blem\la{prop:differentiation}
Let $\cB\subset\Aut(X)$ be 
a Borel subgroup. 
Then $\fb:=\Lie(\cB)$ is a locally integrable 
Borel subalgebra.
\elem
\bproof We have $\cB=\varinjlim \cB_i$, 
where the $\cB_i$ are connected 
algebraic subgroups, see 
Lemma~\ref{lem:nested}.
Any vector $v\in \fb=T_e\cB$ is tangent to 
$\cB_i$ for $i\gg 1$.
Therefore, $\fb=\varinjlim \fb_i$, 
where the $\fb_i=\Lie(\cB_i)$
are locally finite, 
see \cite[Section~1.6, Corollary]{KZ24}.  

Let $\fh=\varinjlim \fh_i\subset \Lie(\Aut(X))$ be 
a solvable locally integrable subalgebra
 that contains~$\fb$, where 
the $\fh_i$ are locally finite. 
 Let also $G=G_{\min}(\fh) \subset \Aut(X)$ be 
the unique smallest solvable connected 
nested ind-subgroup
such that $\fh \subset \Lie(G)$,
see Proposition \ref{prop:solv}. 
We have 
$G=\varinjlim G_i$, where 
$G_i=G_{\min}(\fh_i)\subset\Aut(X)$ 
is the minimal connected algebraic subgroup
containing $\fh_i$, see Lemma \ref{lem:new}. 
For every $i$, there exists $j$
such that
\[\fb_i=\Lie(\cB_i)\subset \fh_j
\subset \Lie(G_{\min}(\fh_j)).\] 
An embedding of algebraic 
Lie subalgebras implies
an embedding of the corresponding 
connected algebraic subgroups, see 
\cite[Theorem~13.1]{Hum81}.
Therefore, 
$\fB_i\subset G_j=G_{\min}(\fh_j)$, 
and so $\fB\subset G$.
Since  $G$ is solvable, 
$G=\cB$ by the maximality of $\cB$. 
Thus, $\fb=\fh$, which proves the claimed 
maximality of $\fb$.
\eproof
The main result of this section 
is the following theorem. 
\bthm\la{thm:integration}
Let 
$\fb=\varinjlim \fb_i\subset \Lie(\Aut(X))$ 
be a locally integrable Borel subalgebra.  
Then there exists a unique 
Borel subgroup $\cB\subset\Aut(X)$
such that $\fb=\Lie(\cB)$. Moreover, any 
solvable closed connected subgroup 
$H\subset\Aut(X)$
with $\Lie(H)=\fb$ coincides with $\cB$. 
\ethm
\brem There exists an example of a
closed ind-subgroup 
$G\subset\Aut(\AA^3)$ 
and a proper closed subgroup 
$H\subset G$ 
such that $\Lie(H)=\Lie(G)$,
see \cite[Theorem~17.3.1]{FK18}.
Therefore, the uniqueness in Theorem 
\ref{thm:integration} is a 
special phenomenon. 
\erem
We  prove the existence 
and the uniqueness statements 
separately.
\bproof[Existence proof]
 By Lemma \ref{lem:new},
there exists a unique minimal 
solvable connected algebraic subgroup 
$G_i\subset\Aut(X)$ such that
$\fb_i\subset \Lie(G_i)$. 
We have $G_i=T_i\ltimes U_i$,
where $T_i$ is a maximal torus and $U_i$
 is the unipotent radical of $G_i$. 

By the uniqueness of the subgroups $G_i$ 
we have $G_i\subset G_{i+1}$
and $U_i \subset U_{i+1}$. Therefore,
$U:=\varinjlim U_i$ is 
a nested unipotent ind-subgroup of $\Aut(X)$.
The maximal tori $T_i$ of $G_i$ 
can be chosen such that 
$T_i\subset T_{i+1}$.
Since the sequence $(T_i)_i$ 
stabilizes, we can assume that 
$T_i=T_1=:T$ for all $i$. Thus, 
$G:=\varinjlim G_i=T\ltimes U\subset\Aut(X)$ 
is a connected nested ind-subgroup. 
By Theorem \ref{thm:closed}, $G$ is closed. 

Since
$U$ is solvable, see  \cite[Theorem~5.3.1]{KZ24},
$G=T\ltimes U$ is also solvable. 
We have 
\[\Lie(G)=\varinjlim \Lie(G_i)=\Lie(T)\ltimes 
\varinjlim \Lie(U_i)=\Lie(T)\ltimes \Lie(U),\] 
where the factors are solvable,
see Theorem \ref{thm:nested}. 
Therefore, $\Lie(G)$ is 
a locally integrable solvable 
subalgebra of 
$\Lie(\Aut(X))$ filtered by 
the locally finite subalgebras $\Lie(G_i)$. 
By assumption,
$\fb=\varinjlim \fb_i\subset\Lie(G)$ is 
a locally integrable Borel subalgebra of 
$\Lie(\Aut(X))$.
Therefore, we have $\fb=\Lie(G)$. 

Let $\hat G\subset\Aut(X)$ be 
a solvable closed connected subgroup 
 containing $G$. By Theorem \ref{thm:nested}, 
 $\hat G$ is nested,
and so $\Lie(\hat G)\subset\Lie(\Aut(X))$ is 
a solvable locally integrable Lie subalgebra 
 containing $\Lie(G)=\fb$. 
 Since $\fb$
 is a locally integrable Borel subalgebra  
 of $\Lie(\Aut(X))$, 
 we have $\Lie(\hat G)=\fb$. 
 
 Let $\{\hat G_\alpha\}_{\alpha\in\cA}$ 
 be a totally ordered ascending chain of 
 solvable closed connected subgroups
 containing $G$. 
 By Theorem \ref{thm:nested}, 
 they are all nested by connected 
 algebraic subgroups and have
derived lengths $\le\dim(X)+1$. Consider
the closure $\bar S$ of the solvable subgroup
$S=\cup_{\alpha\in\cA} \hat G_\alpha$, 
see Lemma \ref{lem:solv-union}. 
Then $\bar S$ is
solvable (see \cite[Lemma~5.1.3(2)]{KZ24}) 
closed connected subgroup 
of derived length $\le\dim(X)+1$ with 
$\Lie(\bar S)=\fb$ that majorates our chain. 
 By Zorn's Lemma, 
 the set of all solvable closed connected 
 subgroups 
 of $\Aut(X)$ 
 containing $G$ has a maximal element, 
 say $\cB$
 (cf.  \cite[Proposition~3.10]{FP18}). 
Thus,  $\cB\supset G$ is a solvable closed connected 
subgroup of $\Aut(X)$  
maximal among the subgroups
 of $\Aut(X)$ 
with the same properties. That is, 
$\cB$ is a Borel subgroup of $\Aut(X)$
 such that $\Lie(\cB)=\fb$. 
 \eproof
 \bproof[Uniqueness proof]
Let $\cB$ be a Borel subgroup
 of $\Aut(X)$ with $\Lie(\cB)=\fb$, 
 and let $H\subset\Aut(X)$ be 
 a solvable closed connected subgroup
with $\Lie(H)=\fb$. 
 By Lemma \ref{lem:nested} resp. by 
 Theorem \ref{thm:nested}, 
$\cB=\varinjlim \cB_i$ 
resp. $H=\varinjlim H_i$ 
is nested by solvable connected 
algebraic subgroups $\cB_i$ resp. $H_i$.
This produces the ascending filtrations
$\fb=\varinjlim \fb_i=\varinjlim \fh_i$ 
by 
locally finite Lie subalgebras 
$\fb_i:=\Lie(\cB_i)$ resp.
$\fh_i:=\Lie(H_i)$, 
see Remark \ref{rem:loc-fin}. 
For any $i$ there 
exists $j$ such that 
 $\fh_i\subset \fb_j$ and  
 $\fb_i\subset\fh_{j}$. 
 Therefore,
we have $H_i\subset \cB_j$ and 
$\cB_i\subset H_j$. 
It follows that
 $H=\cB$. 
 \eproof
The following immediate corollary of 
Theorem \ref{thm:integration} 
and 
Lemma \ref{prop:differentiation} 
proves Theorem \ref{thm:1.2}.
\bcor\la{cor:corresp}
The conjugacy classes of the 
Borel subgroups 
of $\Aut(X)$ are in natural
one to one correspondence with 
the $\Ad$-conjugacy classes 
of locally integrable Borel subalgebras
of $\Lie(\Aut(X))$. 
\ecor
\brem Let $\cB\subset\Aut(X)$  
be a Borel subgroup. 
It is not true, in general, that
 $\Lie(\cB)$ is a Borel subalgebra 
of $\Lie(\Aut(X))$, that is, 
a maximal solvable subalgebra 
of $\Lie(\Aut(X))$. 
See Proposition \ref{mthm-triang} 
below
for the corresponding example. 
\erem
\section{Borel subalgebras of 
$\Lie(\Aut(\AA^n))$  for small $n$}\la{sec:4}
In this section we prove Theorem \ref{thm:1.1}, 
see Propositions \ref{mthm-1} ($n=2$) 
and \ref{mthm-triang} ($n=3$).
\subsection{Borel subgroups of $\Aut(\AA^n)$} 
Recall that the lower and upper triangular 
(de Jonqu\`eres) subgroups 
of $\Aut(\AA^n)$ 
are
\[{\JONQ}^+(\AA^n)=\{\phi\colon (x_1,\ldots,x_n)
 \longmapsto (\alpha_1 x_1+p_1,\ldots,
 \alpha_n x_n+p_n)\,|\,\alpha_i
\in\kk^*, \,
p_i\in\kk[x_{i+1},\ldots,x_n]\},\]
resp.
\[{\JONQ}^-(\AA^n)=\{\phi\colon (x_1,\ldots,x_n)
 \longmapsto (\alpha_1 x_1+p_1,\ldots,
 \alpha_n x_n+p_n)\,|\,\alpha_i
\in\kk^*, \,
p_i\in\kk[x_{1},\ldots,x_{i-1}]\}.\]
They are conjugate via the reversion 
$(x_1,\ldots,x_n)\mapsto (x_n,\ldots,x_1)$.

The following results are well known: see 
\cite[Theorem~1.1~and~Corollary~1.2]{FP18} 
for (a), (d) and (e), 
\cite[Corollary~1.5]{RUVS25} for (b), and 
\cite[Theorem~1]{BEE16} for (c).
\bthm\la{thm:B-sbgrps} $\,$

\bnum
\item[$(a)$] For all $n \ge 2$, 
the triangular de Jonqui\`eres subgroups 
$\JONQ^\pm(\AA^n)$ are Borel subgroups 
maximal among the solvable subgroups 
of $\Aut(\AA^n)$.

\item[$(b)$] For all $n \ge 3$, 
there exist Borel subgroups 
of $\Aut(\AA^n)$ non-conjugate to
$\JONQ^\pm(\AA^n)$.

\item[$(c)$] Every Borel subgroup of 
$\Aut(\AA^2)$ is conjugate to both 
$\JONQ^\pm(\AA^2)$. 

\item[$(d)$] The subgroup 
$\JONQ^\pm(\AA^2)$ 
is maximal among the closed subgroups 
of $\Aut(\AA^2)$.

\item[$(e)$] The subgroup 
$\JONQ^\pm(\AA^2)$ 
is not maximal among all proper subgroups of 
$\Aut(\AA^2)$.
\enum
\ethm
Note that the locally integrable Borel 
subalgebra 
\[\fj^\pm_n:=\Lie({\JONQ}^\pm(\AA^n))
\subset\Der(\kk[x_1,\ldots,x_n])\]
consists of the lower resp. upper  triangular 
derivations.

\subsection{The case of the affine line}\la{ss:4.1}
Here we classify the Borel subalgebras of 
the Lie algebra $\Der(\kk[x])$.
The polynomials in $\kk[x]$ of the form
\[f(x)=\lambda (x-\alpha)^k+\mu (x-\alpha),
\quad\text{where}\quad k\ge 0
 \quad\text{and}\quad
\alpha,\lambda,\mu\in\kk,\]
will be called \emph{special}.
Non-special polynomials exist 
in any degree $\ge 3$. 
\bprop\la{mthm1} $\,$

\bnum
\item[$(a)$] Every Borel subalgebra 
of  $\Der(\kk[x])$ 
is of dimension $1$ or $2$.
\item[$(b)$] 
The one-dimensional Borel subalgebras 
of  $\Der(\kk[x])$ are precisely
the subalgebras $\kk p(x)\p/\p x$, 
where $p\in\kk[x]$ is non-special.
\item[$(c)$] The two-dimensional 
Borel subalgebras 
of  $\Der(\kk[x])$ are precisely 
the metabelian  Lie subalgebras
 \be\la{eq:dim-2}
\left(\kk (x-\alpha)^k\oplus\kk 
 (x-\alpha)\right)
 \p/\p x,\ee
where $\alpha\in\kk$, 
$k\ge 0$ and $k \neq 1$.
\enum
\eprop
See the proof at the end of this subsection. 
\bcor Any derivation 
$f(x)\p/\p x\in\Der(\kk[x])$ 
with a non-special polynomial $f\in\kk[x]$ 
is contained in the unique
Borel subalgebra $\kk f(x)\p/\p x$ of 
$\Der(\kk[x])$.
\ecor
\brem\la{rem:4.2} 
Letting in \eqref{eq:dim-2}  
$\alpha=0$ and $k=0$ 
we see that
\[\Lie(\Aut(\AA^1))=
\kk x\p/\p x\ltimes \kk \p/\p x\]
 is a Borel subalgebra of $\Der(\kk[x])$.
 In fact, it is a unique locally integrable Borel 
 subalgebra of $\Der(\kk[x])$.
\erem
Recall that
\[[f \p/ \p x, g  \p/ \p x] = (fg'-f'g)  \p/ \p x.\]
The proof of Proposition \ref{mthm1} 
is preceded by the following notation 
and Lemma \ref{lem-orders}.  
\bnota 
Given $f\in\kk[x]$ we set $\nu(f)$ 
stand for the zero order  at $x=0$ of $f$. 
If $f, g\in\kk[x]$ and 
$\nu(f)\neq\nu(g)$, 
then
\be\la{eq1} \nu(fg'-f'g)=\nu(f)+\nu(g)-1.\ee 
Especially,
if $\nu(f), \nu(g)\ge 2$, then 
\[\nu(fg'-f'g)>\max\{\nu(f),\nu(g)\}.\]
\enota
\blem\la{lem-orders}
Consider a subalgebra 
$\fg\subset\Der(\kk[x])$.
Let $u_i=f_i(x)\p/\p x\in\fg$, $i=1,2$ be such that 
$\min\{\nu(f_1),\nu(f_2)\}\ge 2$ and 
$\nu(f_1)\neq\nu(f_2)$. 
Then the commutant
 $[\fg,\fg]$ contains  two elements 
 $v_i=g_i(x)\p/\p x$, $i=1,2$ such that 
$\min\{\nu(g_1),\nu(g_2)\}\ge 3$ and 
$\nu(g_1)\neq\nu(g_2)$. 
\elem
\bproof Consider the elements 
\[v_1=[u_1,[u_1,u_2]]=g_1(x)\p/\p x\quad
\text{and}\quad v_2=[u_1,u_2]=g_2(x)\p/\p x.\] 
Using \eqref{eq1}, we obtain:
\[\nu(g_1)=2\nu(f_1)+\nu(f_2)-2>\nu(g_2)
=\nu(f_1)+\nu(f_2)-1 \ge 3.\]
\eproof
\bproof[Proof of Proposition \ref{mthm1}] 
Let $\fg\subset\Der(\kk[x])$ be a
Lie subalgebra of dimension $\ge 4$, 
and let  $u_i=f_i(x)\p/\p x\in\fg$, $i=1,2,3,4$ 
be linearly independent. Suppose that 
\[\nu(f_i)\le \nu(f_{i+1}),\quad
i=1,2,3.\] 
By replacing the $u_i$ by their 
appropriate linear combinations, 
we may assume that 
\[\nu(f_i)<\nu(f_{i+1}),\quad
i=1,2,3.\] 
It follows that
\[2\le \nu(f_3)<\nu(f_4).\]
Applying Lemma \ref{lem-orders}, 
we see that all successive commutants 
of $\fg$ 
do not vanish, so $\fg$ is not solvable. 
Thus, every solvable subalgebra, 
and in particular 
every Borel subalgebra of $\Der(\kk[x])$, 
is of dimension $\le 3$.
Finite dimensional Lie subalgebras 
of $\Der(\kk[x])$ 
are classified in \cite[Proposition~3]{AMP11}. 
Every one-dimensional Lie subalgebra  
is abelian.
Up to isomorphism, there is a unique 
two-dimensional Lie subalgebra; 
it is not abelian but
metabelian, and therefore solvable. 
Up to isomorphism, there exists a unique 
subalgebra $\fh$ of dimension 3,
 isomorphic to $\fs\fl(2,\kk)$, 
see  Remark \ref{rem:sl2}.  This proves (a).
Now, (b) and (c) 
follow from the explicit description 
of these subalgebras,
see \cite[Proposition~3]{AMP11}. 
\eproof
\subsection{The case of the affine plane}
\la{ss:4.2}
\subsubsection{Borel subgroups}
 \la{sss:Borel-sbgrp}
We consider the following remarkable subgroups of 
$\Aut(\AA^2)$ (resp. of
 $\Aut_0(\AA^2)$):
 \begin{itemize}
 \item the subgroup $\mathfrak{U}_2\subset\Aut(\AA^2)$ 
 of automorphisms of $\AA^2$ 
 whose Jacobian 
is equal to 1;
 \item the triangular 
de Jonqui\`eres subgroups
 \[{\JONQ}^+(\AA^2)=\{\phi\colon (x,y)
 \longmapsto (\alpha
x+p(y),\beta y+\gamma)\,|\,\alpha,
\beta\in\kk^*, \gamma\in\kk,\,
p\in\kk[t]\}\]
and
 \[{\JONQ}^-(\AA^2)=\{\phi\colon (x,y)
 \longmapsto (\alpha
x+\gamma,\beta y+q(x))\,|\,\alpha,
\beta\in\kk^*, \gamma\in\kk,\,
q\in\kk[t]\}\]
\[\text{(resp.}\quad  {\Jonq}^\pm(\AA^2)=
\{\phi\in{\JONQ}^\pm(\AA^2)\,|\,\phi(0)=0\}),\]
\item the standard maximal torus
\[\TT=\{\psi\colon(x,y)\longmapsto (\alpha
x, \beta y)\}\subset 
{\Jonq}^+(\AA^2)\cap{\Jonq}^-(\AA^2);\]
\item the unipotent radical
\[U^\pm=\{\phi\in{\JONQ}^\pm(\AA^2)\,|\,
\alpha=\beta=1\}
\quad\text{(resp.}\quad
U^\pm_0=\{\phi\in{\Jonq}^\pm(\AA^2)
\,|\,\alpha=\beta=1\}).\]
\end{itemize}
We have
\[{\JONQ}^\pm(\AA^2)=\TT\ltimes 
U^\pm\quad(\text{resp.}
\quad {\Jonq}^\pm(\AA^2)=
\TT\ltimes U^\pm_0).\]
\subsubsection{$\Lie(\Aut(\AA^2))$
 as a bigraded Lie algebra}

Recall that $\Lie(\Aut(\AA^2))=
\Vect^{\mathrm c}(\AA^2)$ 
consists of polynomial vector fields 
on $\AA^2$ 
with constant divergence, 
where the divergence of a vector field 
$\delta=f \p/\p x + g\p/\p y$ is 
${\rm div}(\delta)= \p f/\p x+\p g/\p y$.
The Lie subalgebra $\Vect^0(\AA^2)\subset
 \Vect^{\mathrm c}(\AA^2)$ of vector fields
over $\AA^2$ with zero divergence coincides 
with $\Lie(\mathfrak{U}_2)$, 
see for example \cite[Proposition~15.7.2]{FK18}.
There is a decomposition
\be\la{eq:vec} {\Vect}^{\mathrm c}(\AA^2)=
{\Vect}^0(\AA^2)\oplus 
\kk(x \p/\p x + y \p/\p y).\ee
Letting 
\[{\Vect}_0^0(\AA^2):=
{\Vect}^0(\AA^2)\cap 
{\Vect}_0^{\mathrm c}(\AA^2),
\quad\text{where}\quad
{\Vect}_0^{\mathrm c}(\AA^2)=
\{\delta\in {\Vect}^{\mathrm c}(\AA^2)
\,|\,\delta(0)=0\},\]
we obtain
\[\Lie({\Aut}_0(\AA^2))=
{\Vect}_0^0(\AA^2)\oplus 
\kk(x \p/\p x + y \p/\p y),\]
see \cite[Section~6.2]{KZ24}.

The standard $\ZZ^2$-grading 
on $\kk[x,y]$ induces 
a $\ZZ^2$-grading on 
$\Lie(\Aut(\AA^2))=\Vect^{\mathrm c}(\AA^2)$
and on $\Lie(\mathfrak{U}_2)=\Vect^0(\AA^2)$.
 For $i,j\ge 0$ we have
\[\bideg (x^iy^j\p/\p x)=(i-1,j)\quad\text{and}
\quad\bideg (x^iy^j\p/\p y)=(i,j-1).\]
Thus, the bidegrees of homogeneous derivations 
are the pairs $(m,n)\in\ZZ^2$, 
where $m,n\ge -1$ and $(m,n)\neq (-1,-1)$. 
If $\p\in \Der (\kk[x,y])$ is a homogeneous 
derivation of bidegree 
$(a,b)$ and $x^iy^j\notin \ker(\p)$, then 
$\bideg(\p(x^iy^j))=(a+i,b+j)$. 
For non-commuting homogeneous derivations 
$\p$ and $\delta$ we have
$\bideg([\delta,\p])=\bideg(\delta)+\bideg(\p)$.

The homogeneous elements 
of $\Vect^0(\AA^2)$ with bidegree $(a,b)$, 
where $ a,b\ge 0$, 
are (up to a scaling) of the form 
\[\p_{a,b}
=x^ay^b((b+1)x\p/\p x- (a+1)y\p/\p y).\]
For $\min\{a,b\}=-1$ they are of the form
\[\p_{-1,n}=(n+1)y^n\p/\p x
\quad\text{and}\quad 
\p_{m,-1}=-(m+1)x^m\p/\p y.\]
Every homogeneous 
locally nilpotent derivation $\p$
of $\kk[x,y]$ belongs to ${\Vect}^0(\AA^2)$. 
Up to a scalar factor, $\p$ coincides 
with $\p_{-1,n}$ or $\p_{m,-1}$.
These homogeneous derivations 
verify the following
commutation relations: 
\be\la{eq:comm}
[\p_{a,b},\p_{a',b'}]=
\det\left(\begin{matrix} 
a'+1&a+1\\b'+1&b+1\end{matrix}\right)
\p_{a+a', b+b'}\quad\text{for all} \,\,\,
a,b,a',b'\ge -1, \,\,\, (a,b),\,(a',b') \neq (-1,-1),\ee
where we let $\p_{-1,-1}=0$. 
Every semisimple homogeneous derivation
of ${\rm Vec}^{\mathrm c}(\AA^2)$
has bidegree $(0,0)$. 
They are of the form 
$\delta_{\alpha,\beta}=\alpha x\p/\p x
+\beta  y\p/\p y$ and satisfy 
the commutation relations 
\be\la{eq:comrel}
[\delta_{\alpha, \beta}, 
\delta_{\gamma, \eta}] 
= 0\quad\text{and}\quad
[\delta_{\alpha, \beta}, \p_{a,b}]= 
(\alpha a + \beta b) \p_{a, b}.\ee
Furthermore, we have 
\be\la{eq:concord} 
\p_{a,b}(x^ky^l)=\det\left(\begin{matrix} 
k&a+1\\l&b+1\end{matrix}\right)
x^{k+a}y^{l+b}.\ee
\subsubsection{Borel subalgebras of 
$\Lie(\Aut(\AA^2))$}
Let us begin with the following remark. 
\brem\label{rem:Lie-A2} 
Let us show that $\Lie(\Aut(\AA^2))$ 
contains non-conjugate Borel subalgebras. 
Let $\fg\subset\Lie(\Aut(\AA^2))$ 
be a Borel subalgebra  
$\Ad$-conjugate to the triangular subalgebra 
$\fj^+_2:=\Lie(\JONQ^+(\AA^2))$. 
Then $\fg$ is locally integrable,  
and therefore consists of 
 locally finite derivations.
However, 
$\Lie(\Aut(\AA^2))=
\Vect^{\mathrm c}(\AA^2)$ 
contains elements that 
are not locally finite;
for example, $\p=\p_{1,0}
=x^2\p /\p x - 2xy \p / \p y\in\Vect^{0}(\AA^2)$ 
is such an element. 
Indeed, the degrees 
of the iterates $\p^n(x)$ are unbounded. 
By Corollary \ref{cor:maximal},
there exists a Borel subalgebra 
$\ff\subset\Lie(\Aut(\AA^2))$
containing $\p$. Therefore, $\ff$ is
not conjugate to 
$\fj^+_2$ and does not correspond to 
a Borel subgroup of $\Aut(\AA^2)$.
\erem
Consider the following Lie subalgebras 
of $\Lie(\Aut(\AA^2))$ 
(resp. $\Lie(\Aut_0(\AA^2))$):
  \begin{itemize}
 \item 
the metabelian (resp. abelian) 
Lie subalgebras 
of triangular locally nilpotent derivations 
\[\fu^+_2=\Lie(U^+)= \kk[y]\p / \p x \oplus
\,\kk\p/\p y\quad\text{and}\quad 
\fu^-_2=\Lie(U^-)=\kk\p / \p x\oplus
\,\kk[x]\p/\p y\]
\[\text{(resp.}\quad
\fu^+_{2,0}=\Lie(U^+_0)=y\kk[y] \p / \p x
\quad\text{and}\quad 
\fu^-_{2,0}=\Lie(U^-_0)=x\kk[x] \p / \p y);\] 
\item the abelian toral Lie subalgebra
\[\ft_2=\Lie(\TT)=\kk x\p/\p x\oplus
\kk y\p/\p y;\]
\item the triangular Lie subalgebras 
of derived length 3 (resp. 2)
\[\fj_2^\pm=\Lie({\JONQ}^\pm(\AA^2))
=\Lie(\TT\ltimes U^\pm)=
\ft_2\ltimes \fu_2^\pm\]
\[\text{(resp.}
\quad \fj_{2,0}^\pm=
\Lie({\Jonq}^\pm(\AA^2))
=\Lie(\TT\ltimes U^\pm_0)=
\ft_2\ltimes \fu^\pm_{2,0}),\]
where
\be\la{eq:decomp-j2} \fj_2^+=\kk[y]\p / \p x\oplus 
\kk x\p / \p x\oplus \kk y\p / \p y
\oplus\kk \p / \p y\ee
\be\la{eq:decomp-j2-0} 
\text{(resp.}\quad \fj_{2,0}^+
=y\kk[y]\p / \p x
\oplus \kk x\p / \p x\oplus \kk y\p / \p y).\ee
\end{itemize}
Using 
\eqref{eq:comm}
 we see that
\[\Lie(\mathfrak{U}_2)={\Vect}^0(\AA^2)
=\langle \fu^+_2,\,\fu^-_2\rangle_{\Lie}.\]
Therefore, by \eqref{eq:vec} we have
\[\Lie(\Aut(\AA^2))=
{\Vect}^{\mathrm c}(\AA^2)
=\langle \fu^+_2,\,\fu^-_2\rangle_{\Lie}
\oplus\kk(x\p/\p x+y\p/\p y).\]
Similar equalities hold for all 
$\AA^n$ with $n\ge 2$; 
see, for example, the proof 
of \cite[Lemma~2]{Sha81}.

In the proof of Proposition \ref{mthm-1} below
we use the following simple lemma.
\blem\la{lem-simple}
Let $A$ be a diagonalizable
endomorphism of a vector space $V$.
Consider the decomposition 
$V=\oplus_i V_i$, where the $V_i$ are
eigenspaces of $A$ 
corresponding to different eigenvalues
$\lambda_i$, $i=1,\ldots,m$.  
Let a subspace  
$U\subset V$ be invariant under $A$, and let
$u=u_1+\cdots+u_m\in U$, where $u_i\in V_i$. 
Then $u_i\in U$ for all $i$. 
\elem
\bproof By symmetry,
it suffices to show that $u_1\in U$. 
The latter is obviously true if $u=u_1$.
Arguing by recursion, we can assume that 
our claim is true for $u=u_1+\cdots+u_{k-1}$. 
Now suppose that $u=u_1+\cdots+u_k$,
where $1<k\le m$.
We have 
\[U\ni Au-\lambda_k u
=\sum_{j=1}^{k-1} (\lambda_j- \lambda_k)u_j.\] 
By the inductive 
 hypothesis, $ (\lambda_j- \lambda_k)u_j\in U$
 for $j=1,\ldots,k-1$, where 
 $\lambda_j- \lambda_k\neq 0$.
  In particular, $u_1\in U$.
\eproof
 \bprop\la{mthm-1} The subalgebra 
 $\fj^\pm_2~=~\Lie(\JONQ^\pm(\AA^2))$ 
 is a Borel subalgebra of $\Lie(\Aut(\AA^2))$.
\eprop
 \bproof 
 To prove that $\fj^+_2$
 is a Borel subalgebra, 
 it suffices to show that
 any Lie subalgebra
 $\fh\subset\Lie(\Aut(\AA^2))$ 
properly containing $\fj^+_2$ 
 also contains a subalgebra
 isomorphic to $\fs\fl(2,\kk)$.
Therefore, $\fh$ cannot be solvable. 
 
Choose 
 $v\in\fh\setminus\fj^+_2$. 
As a vector space, $\fj^+_2$ is generated 
by its homogeneous elements.
Hence, $v$ contains homogeneous 
components that do not belong to $\fj^+_2$.
We can assume that none of 
the homogeneous 
components of $v$ belong to $\fj^+_2$. 

Since $\ft_2\subset\fj^+_2$, 
the semisimple derivations 
$A=\ad(x\p/\p x)$ and $B= \ad(y\p/\p y)$ 
from $\End(\Der(\kk[x,y]))$ leave $\fh$
 invariant. 
Let $v_{a,b}$ 
be the homogeneous component of $v$
of bidegree $(a,b)$. 
Then $v_{a,b}$ is an eigenvector
of $A$
(resp. of $B$)
with eigenvalue $a$ (resp. $b$). 
Let us show that $v_{a,b}\in\fh\setminus\fj^+_2$.

Note that 
$A|_{\fh}$ is locally finite. Therefore, $v$ 
belongs to a finite-dimensional 
$A$-invariant vector subspace 
$V\subset\fh$. Being semisimple, 
$A|_V$ is diagonalizable and
defines an eigenspace
decomposition of $V$. 
If $v=\sum_{(a,b)} v_{a,b} $ 
is the usual weight decomposition,
then the $A$-homogeneous 
components of $v$ are of the form 
 $v_a=\sum_b v_{a,b}$. 
They belong to $V$ due to 
Lemma \ref{lem-simple}.
The same argument applied to 
$B=\ad(y\p/\p y)\in\End(\Der(\kk[x,y]))$
shows that the $B$-homogeneous 
components $v_{a,b}$ of $v_a$ 
also belong to $V$. 
Thus, $v_{a,b}\in \fh\setminus \fj^+_2$
for all $(a,b)$. 

Let us return to the proof of the proposition.
First, note  that $v_{a,b}\in\kk^*\p_{a,b}$.
Consequently, $\p_{a,b}\in \fh\setminus\fj^+_2$ for some
pair $(a,b)\in\ZZ^2$ such that $a\ge 0$, $b\ge -1$, 
and $(a,b)\neq (0,0),\, (0,-1)$. Thus,
it suffices to examine the following possibilities 
for the bidegree $(a,b)$ of $\p_{a,b}\in\fh$:
\begin{itemize}
\item $(a,b)=(a,-1)$, with $a\ge 1$. 
Since $\p_{-1,0}=\p/\p x\in\fj^+_2$,
formula \eqref{eq:comm} implies
\[\p_{a-1,-1}\in\kk^*[\p_{-1,0}, \p_{a,-1}]
\subset\fh.\]
By iteration, we deduce that $x\p/\p y\in\fh$. 
Since $y\p/\p x\in\fj^+_2\subset\fh$, we have
\[\fs\fl(2,\kk)\simeq\langle x\p/\p y,
 \,y\p/\p x\rangle_{\Lie}\subset\fh.\]
\item $a\ge 1$, $b\ge 0$. Since $\p_{0,-1}=\p/\p y\in\fj^+_2$, 
we obtain by  \eqref{eq:comm},
\[\p_{a,b-1}\in\kk^*[\p_{0,-1},\p_{a,b}]\subset\fh.\]
Iterating, we arrive at
the case where $b=-1$. 
So, $\p_{a,-1}\in\fh$ for some $a\ge 1$.
This brings us back to the previous case.
\item $(a,b)=(0,b)$, where $b\ge 1$. 
Then, by  \eqref{eq:comm},
\[\p_{0,b-1}\in\kk^*[\p_{0,-1}, \p_{0,b}]\subset\fh.\]
Iterating, we arrive at the case 
where $b=1$, that is,
$\p_{0,1}\in\fh$. 
We then have
\[\fs\fl(2,\kk)\simeq \kk\p_{0,-1}
\oplus \kk \p_{0,0} 
\oplus \kk \p_{0,1}
\subset\fh.\]
 \end{itemize}
 This proves the assertion for 
 $\Lie(\JONQ^+(\AA^2))$;
 the proof  is similar for $\fj^-_2$.
\eproof
\brem Let us show that
the triangular subalgebra 
$\fj_2^+\subset
\Lie(\Aut(\AA^2))$
is not a Borel 
subalgebra of $\Der(\kk[x,y])$. 
It suffices to verify that the subalgebra
\[\fg:=\langle \fj^+_2, \delta
\rangle_{\Lie}
\subset\Der(\kk[x,y]),\quad
\text{where}\quad \delta=xy\p/\p x,\]
 strictly containing $ \fj^+_2$
 is solvable. It is easy to see
that 
\[[\fg,\fg] \subset  \fu^+_2 \oplus \kk\delta 
\oplus \kk x \p/\p x
\quad\text{and}
\quad \fg^{(2)}\subset \fj_2^+.\]
Now the claim follows. 
\erem
As per the following proposition,
an analogue of Proposition \ref{mthm-1} 
does not hold for the Lie algebra
$\Lie(\Aut_0(\AA^2))$.
\bprop \la{prop:non-borel}
The de Jonqui\`eres Lie subalgebra
$\fj_{2,0}^+\subset\Lie(\Aut_0(\AA^2))$
 is not a Borel subalgebra of 
$\Lie(\Aut_0(\AA^2))$. More precisely, 
the distinct 
Borel subalgebras of 
$\Lie(\Aut_0(\AA^2))$ that
strictly contain $\fj_{2,0}^+$
form an infinite sequence 
\[\fh_k:=\langle\fj_{2,0}^+,\,\p_{0,k}
 \rangle_{\Lie},\quad k=1,2,\ldots\]
None of the $\fh_k$ are locally integrable.
\eprop
\bproof We have $\p_{0,k}\in\Lie(\Aut_0(\AA^2))$,
and consequently $\fh_k\subset\Lie(\Aut_0(\AA^2))$, 
see \eqref{eq:decomp-j2-0}. 
For a fixed $k\ge 1$, from formula \eqref{eq:concord}
we obtain
\[\p^{(d)}_{0,k}(x)=c_dxy^{kd} 
\quad\text{where}\quad c_d\neq 0,\quad d> 0.\]
It follows that $\p_{0,k}$ is not locally finite. 
This proves the second assertion. 

To show that $\fh_k$ is solvable, note that
$\fii_k :=\fu^+_{2,0}\oplus\kk\,\p_{0,k}$
is an ideal of $\fh_k$ and that $[\fh_k, \fh_k] =\fii_k$. 
Furthermore, $\fu^+_{2,0}$ is an abelian ideal in $\fii_k$,
whence $\fh_k^{(2)}\subset \fu^+_{2,0}$, 
and $\fh_k^{(3)}=0$.

Let us show that each $\fh_k$ is a 
maximal solvable subalgebra of 
$\Lie(\Aut_0(\AA^2))$.
Suppose that a Lie subalgebra 
$\fh\subset\Lie(\Aut_0(\AA^2))$
strictly contains $\fh_k$. 
Proceeding as in the proof 
of Proposition \ref{mthm-1}, we see that
$\fh$ must contain a homogeneous derivation
$\p_{a,b}$ with $(a,b)\neq (0,k)$, 
where $a\ge 0$. If $a=0$, 
then we let $\p_{a,b}=\p_{0,l}$, where $l=b\neq k$.
If $a>0$, then 
we can replace $\p_{a,b}$ 
by 
\be\la{eq:m} \p_{0,l}\in\kk^*\ad(\p_{-1,m})^{(a)} (\p_{a,b}),
\quad\text{where}\quad l=l(m)=
am+b\neq k. \ee

We claim that the Lie subalgebra 
\[\fh_{k,l}:=\langle \p_{0,k},\,\p_{0,l} \rangle_{\Lie}
\subset \fh\]
is not solvable. 
Let $P_{k,l}\subset\NN$ be the sub-semigroup
generated by $k$ and $l$. By \eqref{eq:comm}, 
we have
$\p_{0,n}\in\fh_{k,l}$ for all $n\in P_{k,l}$. 
It is easy to see by induction, using \eqref{eq:comm}, 
that, for any fixed  
$d>0$, the derived ideal
$\fh_{k,l}^{(d)}$ contains $\p_{0,n}$ for 
all sufficiently large $n\in P_{k,l}$, 
hence it is nonzero. 

This proves that $\fh_{k}$ is a Borel subalgebra 
of $\Lie(\Aut_0(\AA^2))$ strictly containing 
$\fu^+_{2,0}$. By the preceding, $\fh_{k}$ 
cannot contain
$\p_{0,l}$ for $l\neq k$, hence $\fh_{k}\neq \fh_{l}$.

Finally, let us show that the $\fh_{k}$ 
are the only 
Borel subalgebras of $\Lie(\Aut_0(\AA^2))$ 
containing $\fj^+_{2,0}$. 
Indeed, consider a Lie subalgebra 
$\fh\subset\Lie(\Aut_0(\AA^2))$ 
strictly containing 
$\fu^+_{2,0}$. Suppose,   for 
the sake of contradiction, 
that $\fh$ is solvable and that
$\fh\neq\fh_k$ for all 
$k\ge 1$. As before, $\fh$ admits 
a vector basis consisting 
of homogeneous derivations $\p_{a,b}$, 
and there exists $\p_{a,b}\in\fh$ 
such that $a\ge 1$. 
Then, by \eqref{eq:m}, it follows that 
also $\p_{0,am+b}\in\fh$ 
for all $m\ge 1$, which, 
as we have already noted, 
contradicts
the assumption that $\fh$ is solvable. 
\eproof
\subsection{The case 
of $\AA^3$}\la{ss:4.3}
Unlike the affine plane case, 
the triangular subalgebra of 
$\Lie(\Aut(\AA^3))$
is not maximal among 
the solvable subalgebras 
of $\Lie(\Aut(\AA^3))=
\Vect^{\mathrm c}(\AA^3)$. 
In other words, it is a locally integrable 
Borel subalgebra 
(see Theorem \ref{thm:B-sbgrps}(a)) 
and not a Borel subalgebra.

We use the following notation.
\bnota
The ind-subgroup  
$\JONQ^+(\AA^3)$ and its Lie algebra 
$\fj^+_3$  of 
triangular derivations of $\kk[x,y,z]$
are nested and admit
the semidirect product decompositions
\[{\JONQ}^+(\AA^3)=\TT_3\ltimes U^+_3
\quad\text{resp.}\quad
\fj^+_3=\ft_3 \ltimes \fu^+_3,\] 
where 
$\TT_3\simeq\mathbb{G}_{\mathrm m}^3$ 
is the diagonal $3$-torus, $U^+_3$ is 
the unipotent radical of $\JONQ^+(\AA^3)$, 
 \[
\ft_3=\Lie(\TT_3)=\kk 
x\p/\p x \oplus \kk y\p/\p y 
\oplus \kk z\p/\p z \] 
and
\[\fu^+_3=\Lie(U^+_3)=\kk[y,z]\p/\p x \oplus
 \kk[z]\p/\p y \oplus \kk\p/\p z.
\]
\enota
\bprop\la{mthm-triang}
The Lie algebra $\fj^+_3$ 
of triangular 
derivations of $\kk[x,y,z]$ 
is not maximal among the 
solvable subalgebras of 
$\Lie(\Aut(\AA^3))=
\Vect^{\mathrm c}(\AA^3)$. 
\eprop
\bproof Let
\[
\delta=z(x\p/\p x-y\p/\p y)
\in{\Vect}^{\mathrm c}(\AA^3)
\setminus\fj^+_3\,.
\]
Consider the Lie subalgebra 
$\fh=\langle
\fj^+_3, \delta\rangle_{\Lie}
\subset{\Vect}^{\mathrm c}(\AA^3)$  
which strictly contains $\fj^+_3$.
We claim that $\fh$ is also solvable. 
 Indeed, an easy 
calculation shows
that 
\[\ad(\delta)(\kk[y,z]\p/\p x)\subset 
\kk[y,z]\p/\p x,\quad
 \ad(\delta)(\kk[z]\p/\p y)
\subset \kk[z]\p/\p y, \]
and
\[\ad(\delta)(\kk\p/\p z) 
=\kk(x\p/\p x-y\p/\p y).\]
So, 
\be\la{eq:u3}\ad(\delta)(\fu^+_3)\subset
\fu^+_3\oplus\kk(x\p/\p x-y\p/\p y).\ee
Moreover,
\[\ad(\delta)(x\p/\p x)
=\ad(\delta)(y\p/\p y)
=0\quad\text{and}\quad 
\ad(\delta)(z\p/\p z)= -\delta.\] 
Therefore, $\ad(\delta)(\ft_3)=\kk \delta$.
This shows that $\fh=\fj^+_3\oplus\kk \delta$ 
(as vector spaces).
It follows that 
\[[\fh,\fh]=[\fj^+_3,\fj^+_3]\oplus [\fj^+_3,\delta]
\subset \fu^+_3\oplus\kk(x\p/\p x-y\p/\p y)
\oplus\kk\delta.\]
Since $\ad(\delta)(x\p/\p x-y\p/\p y)=0$, 
using \eqref{eq:u3} we deduce
$\fh^{(2)}\subset\fj^+_3$.
Therefore, the Lie algebra 
$\fh$ is solvable.
\eproof
\brem Note that $\delta$ 
is not locally finite. 
Indeed, 
$\delta: xz^n\mapsto xz^{n+1}$.
\erem
\section{Subalgebras of 
$\Lie(\Aut(X_{d,e}))$}\la{sec:5}
This section is devoted to the proof 
of Theorem \ref{thm:1.3}.
\subsection{Generalities} 

Let us recall the notation 
and some results 
from \cite[Section~4]{AZ13} 
and \cite[Section~7.1]{AZ25}.
\bnota\la{nota:gde} Let $\zeta$ be 
a primitive root of unity 
of order $d>1$ and $e$ be
an integer number  with
 $1\le e<d$ and $\gcd(d,e)=1$.
 We denote by
$G_{d,e}$
the cyclic subgroup of $\Aut(\AA^2)$ 
generated by 
$(x,y)\mapsto (\zeta^e x, \zeta y)$.
We consider the 
normal affine toric surface  
\[X_{d,\,e}=\Spec (k[x,y]^{G_{d,\, e}})\]
with the acting torus 
\[\TT':=\TT/G_{d,\, e}
\simeq\mathbb{G}_{\mathrm m}^2.\]
The unique singular point of 
$X_{d,\,e}$ is the image 
of  the origin of $\AA^2$, 
see e.g.~\cite[Section~2.6]{Ful93}. 
Consider the normalizers
\[\cN_{d,\,e}=
{\Norm}_{\Aut(\AA^2)}
(G_{d,\,e})\quad\text{and}
\quad N_{d,e}={\Norm}_{\GL(2,\kk)}
(G_{d,\,e})=\cN_{d,e}\cap\GL(2,\kk).\]
Let $\tau\colon (x,y)\mapsto
 (y,x)$ be the twist. We have
\[N_{d,e}=\begin{cases}
 \GL(2,\kk) & \text{ if   
$e=1$,}\\
\Norm_{\GL(2,\kk)}(\TT)=
\langle\TT,\tau\rangle  &\text{ if  
$ e\neq 1$  and
$e^2\equiv 1\mod d$,}\\
\TT&\text{ if  
$e^2\not\equiv 1\mod d$.}
\end{cases}\]
We consider the lower and 
upper triangular subgroups
${\Jonq}^\pm(\AA^2)$, where
${\Jonq}^-(\AA^2)
=\tau{\Jonq}^+(\AA^2)\tau$, 
see subsection \ref{sss:Borel-sbgrp}.
Let also
 \[ \cN_{d,e}^\pm=\cN_{d,e}\cap 
 {\Jonq}^{\pm}(\AA^2)
 \quad\text{and}\quad 
 N_{d,e}^\pm=N_{d,e}\cap 
 {\Jonq}^{\pm}(\AA^2).\]
 We have
\[\cN_{d,e}^+\cap\cN_{d,e}^-= 
{\Jonq}^{+}(\AA^2) \cap
{\Jonq}^{-}(\AA^2)=\TT\,.\]
Letting $\cB_{d,e}^\pm=
\cN_{d,e}^\pm/G_{d,e}$
we get $\cB_{d,e}^+
\cap \cB_{d,e}^-=\TT'$.
 \enota
\bprop[{\rm see \cite[Lemmas~4.3 
and~4.5]{AZ13}}] 
We have
\be\la{eq:surj} \Aut (X_{d,e}) 
\simeq \cN_{d,e}/G_{d,e},\ee
and $\TT'$
is a maximal torus of $\Aut (X_{d,e})$.
Furthermore, 
\[\cN_{d,e}^+
={\Norm}_{\Jonq^+(\AA^2)}(G_{d,e})
=\{\phi\colon (x,y)\mapsto 
(\alpha x+y^ep(y^d),\beta y)\,|
\,p\in\kk[t],\alpha ,\beta\in\kk^*\,\}\]
and
\[\cN_{d,e}^-
={\Norm}_{\Jonq^-(\AA^2)}(G_{d,e})
=\{\phi\colon (x,y)\mapsto 
(\alpha x,\beta y+x^{e'}q(x^d))\,|
\,q\in\kk[t],\alpha ,\beta\in\kk^*\,\}.\]
The group
\be\la{eq-0}\cN_{d,e}^\pm=\TT
\ltimes U_{d,e}^{\pm}\ee
is metabelian. The commutant
\[[\cN_{d,e}^\pm,\cN_{d,e}^\pm]=
\{\phi\in\cN_{d,e}^\pm\,
|\,\alpha=\beta=1\}\]
is abelian and coincides with 
the unipotent radical 
$U_{d,e}^{\pm}$ of 
$\cN_{d,e}^\pm$.
\eprop
The group $\Aut (X_{d,e})$ 
admits 
a free amalgamated 
product structure.
\bthm[{\rm \cite[Theorem 4.2 
and Proposition 4.4]{AZ13}}]$\,$

\begin{itemize}
\item[$(a)$] 
In case
$e^2\not\equiv 1 \mod d$ 
 we have
\[\cN_{d,e} \simeq \cN_{d,e}^+*_{\TT} 
\cN_{d,e}^-\]
and
\[\Aut (X_{d,e}) \simeq
\cB_{d,e}^+*_{\TT'} \cB_{d,e}^-\,.\]
\item[$(b)$]  In case 
$e^2\equiv 1 \mod d$
we have
\[\cN_{d,e}  \simeq 
{\cN_{d,e}^+}*_{N_{d,e}^+}
N_{d,e}\] 
and
\[\Aut(X_{d,e}) \simeq \cB_{d,e}^+
*_{N_{d,e}^+/G_{d,e}} 
N_{d,e}/G_{d,e}.\]
\end{itemize}
\ethm
\bthm[{\rm \cite[Theorem~1.1]{AZ25}}]
\la{main}  $\,$
\bnum
\item[$(a)$]
Let $e^2\equiv 1\mod d$.
Then the $\cN_{d,e}^{\pm}$  
(resp. the $\cB_{d,e}^{\pm}$)
 are  Borel subgroups 
 of $\cN_{d,e}$ 
 (resp.  of $\Aut(X_{d,e})$), 
 and they are conjugate.
Any Borel subgroup of $\cN_{d,e}$ 
(resp. of $\Aut(X_{d,e})$)
is conjugate to $\cN_{d,e}^+$ 
(resp. to $\cB_{d,e}^+$).
 \item[$(b)$]
Let $e^2\not\equiv 1\mod d$.
Then $\cN_{d,e}$ 
(resp. $\Aut(X_{d,e})$)
has exactly two conjugacy classes of 
Borel subgroups
represented by the  $\cN_{d,e}^{\pm}$
(resp. by the $\cB_{d,e}^{\pm}$).
\enum
 \ethm
\subsection{On the Lie algebra 
$\Lie(\Aut(X_{d,e}))$}
 \blem\la{lem:5.5}
 We have 
 \be\la{eq-1}\Lie(\cN_{d,e}^{\pm})
 =\Lie(\TT)\ltimes\Lie(U_{d,e}^{\pm})
 =\ft_2\ltimes\fu^\pm_{d,e},\ee
where the abelian Lie subalgebras
\be\la{eq-2} \fu^+_{d,e}=\Lie(U_{d,e}^+)
=\{\delta_p 
:= y^ep(y^d) \p / \p x\,|\,p\in\kk[y]\}=
\bigoplus_{k\ge 0} \kk y^{e+kd} \p / \p x\ee
and
\be\la{eq-3}\fu^-_{d,e}=\Lie(U_{d,e}^-)
=\{\delta'_q 
:= x^{e'}q(x^d) \p / \p y\,|\,q\in\kk[x]\}
=\bigoplus_{l\ge 0}
 \kk x^{e'+ld} \p / \p y\ee
 are made up of triangular
  locally nilpotent 
 derivations from $\Lie(\cN_{d,e})$, while
the toral subalgebra 
\be\la{eq-4}\ft_2=\Lie(\TT)=
\langle  x\p / \p x,  
y\p / \p y\rangle_{\Lie}=
\kk x\p / \p x\oplus \kk y\p / \p y\ee
is made up of semisimple 
 derivations from $\Lie(\cN_{d,e})$.
 \elem
 \bproof Let $\fG$ be an ind-group 
 and $G\subset\fG$ be an algebraic 
 subgroup. Then we have 
 ${\Lie(G)\subset\Lie(\fG)}$. 
 If $\fG=\varinjlim G_i$
 is nested by algebraic subgroups 
 $G_i$, then
 ${\Lie(\fG)=\varinjlim \Lie(G_i)}$.
 
 The equality \eqref{eq-1} is 
 a direct consequence of
  \eqref{eq-0}.
 Note that the replica $\delta_p$ 
 (resp. $\delta'_q$)  
 of $\p / \p x$ (resp. of $\p / \p y$)
is a locally nilpotent derivations 
of $\kk[x,y]$. 

 Every $u\in U_{d,e}^+$ 
 (resp. $v\in U_{d,e}^-$) 
 is contained in a $\Ga$-subgroup 
\[ U^+_p :=\exp(\kk\delta_p)=
 \{(x+ ty^ep(y^d), y)\,|\,t\in\kk\}\quad
 \text{(resp.} \quad U_q^-
 :=\exp(\kk\delta'_q)).\]
Now \eqref{eq-2} and 
\eqref{eq-3} follow, 
and \eqref{eq-4} 
is immediate.
 \eproof
 \bcor\la{cor:19} 
 The homogeneous locally 
 nilpotent derivations
\[\p^+_{e+kd} = y^{e+kd}\p / \p x
\quad\text{resp.}
\quad\p^-_{e'+ld}= x^{e'+ld} \p / \p y,
\quad
k,l=0,1,2,\ldots\]
generate the abelian Lie algebra 
$\fu^+_{d,e}$ 
resp. $\fu^-_{d,e}$
as a vector space. 
 \ecor
 \bproof 
 Using \eqref{eq-2} resp. \eqref{eq-3}, 
 for $p(y)=y^k$ 
resp. $q(x)=x^l$ we obtain
 \[ \Lie(U^+_p)=\kk\p^+_{e+kd}\subset
  \fu^+_{d,e}\quad\text{resp.}\quad
 \Lie(U^-_q)=\kk\p^-_{e'+ld}\subset
  \fu^-_{d,e},\quad k,l=0,1,2,\ldots ,\] 
  and therefore
$\fu^+_{d,e}=\bigoplus_{k \ge 0}
  \kk\p^+_{e+kd}$ resp. $
   \fu^-_{d,e} =\bigoplus_{l \ge 0}
   \kk\p^-_{e'+ld}$.  
   \eproof
\brem\la{rem:suggest}
The subgroup $\cN_{d,e}$ is defined 
 as the normalizer of $G_{d,e}$ 
 in $\Aut(\AA^2)$,
 see Notation~\ref{nota:gde}.
However, since the origin of $\AA^2$ 
is the unique fixed point of 
 $G_{d,e}$, we have 
 $\cN_{d,e}\subset\Aut_0(\AA^2)$.
 \erem
Let $G$ be a connected 
algebraic group
and $F\subset G$ be a finite subgroup.
It is well known that the Lie algebra 
of the centralizer $\Cent_G(F)$
coincides with the fixed-point 
subalgebra $\Lie(G)^{\Ad(F)}$, 
see for example
\cite[Chapter~1, Section~2.11, 
Exercise~4]{OV90}. 
In our setup, there is 
the following partial analog, 
cf. \cite[Section~6.2, formula (2)]{KZ24}.
 \blem\la{lem:inclu} 
 We have $\Lie(\cN_{d,e})\subset 
 \Lie(\Aut_0(\AA^2))^{\Ad(G_{d,e})}$.
In other words, the $\Ad$-action of $G_{d,e}$ 
on $\Lie(\Aut(\AA^2))$ 
 is identical on the Lie subalgebra   
 $\Lie(\cN_{d,e})$, and then also on 
 the Lie subalgebras
 ${\ft_2=\Lie(\TT), \,\Lie(U_{d,e}^{\pm})}$, 
 and  $\Lie(\cN_{d,e}^{\pm})$ of $\Lie(\cN_{d,e})$. 
 \elem
\bproof 
We have
\[\Lie(\cN_{d,e})=\Lie(\cN_{d,e}^0)
\quad\text{and}\quad
\cN_{d,e}^0\subset
{\Cent}_{\Aut_0(\AA^2)}(G_{d,e}).\]
The second inclusion follows 
from the fact that the
conjugation map sends $\cN_{d,e}$ 
into the finite group $\Aut(G_{d,e})$;
 it is therefore constant on the connected
component $\cN_{d,e}^0$ 
and coincides with the identity there. 
Now, the inclusion $\Lie(\cN_{d,e})
=\Lie(\cN_{d,e}^0)\subset 
 \Lie(\Aut_0(\AA^2))^{\Ad(G_{d,e})}$ follows. 
\eproof
\brem From the decomposition
\[\Lie(\Aut(\AA^2))=\Lie({\Aut}_0(\AA^2))
\oplus \kk\p/\p x\oplus\kk\p/\p y\]
it follows  that 
$\Lie(\Aut(\AA^2))^{\Ad(G_{d,e})}
=\Lie(\Aut_0(\AA^2))^{\Ad(G_{d,e})}$.
\erem
The following corollary is immediate.
\bcor\la{cor:5.8}
 The surjection 
 $\cN_{d,e} \to \Aut(X_{d,e})
 =\cN_{d,e}/ G_{d,e}$, 
 see \eqref{eq:surj},  
 induces isomorphisms 
\[\Lie(\Aut(X_{d,e}))\simeq
\Lie(\cN_{d,e}/G_{d,e})
\simeq\Lie(\cN_{d,e}).\]
 \ecor
\brems\la{rem:5.9} $\,$

1. The $\TT$-action 
 on $\Vect(\AA^2)$
 leaves invariant the subalgebra 
 $\Lie(\cN_{d,e})$.
By Lemma \ref{lem-simple}, 
for any $\p\in\Lie(\cN_{d,e})$
 every homogeneous component 
 of $\p$ also
 belongs to $\Lie(\cN_{d,e})$.

 2. 
 Since 
 $ \ft_2=\Lie(\TT)\subset\Lie(\cN_{d,e})$,  
 the homogeneous derivations 
 $\alpha x \p/\p x+\beta y \p/\p y\in \ft_2$ 
 and $\p^+_{e+kd}\in\fu_{d,e}^+$, 
 $\p^-_{e'+ld}\in \fu_{d,e}^-$  
 belong to 
 $\Lie(\cN_{d,e})$. 
 \erems
 \blem\la{lem:5.10} 
 A homogeneous derivation 
 $\p\in\Vect^{\mathrm c}(\AA^2)$ 
 of bidegree $(a,b)$
 is fixed by the ${\Ad(G_{d,e})}$-action 
 if and only if $ae+b\equiv 0\mod d$.
 Therefore, we have
 \be\la{eq:long} 
 \Lie({\Aut}_0(\AA^2))^{\Ad(G_{d,e})}=
 \kk x\p/\p x \oplus\kk y\p/\p y\oplus
 \bigoplus_{k\ge 0} \kk\p^+_{e+kd}\oplus
 \bigoplus_{l\ge 0}\kk\p^-_{e'+ld}
\oplus \bigoplus_{(a,b)\in
\Lambda_{d,e}} \kk\p_{a,b},\ee
 where the submonoid
 $\Lambda_{d,e}\subset \ZZ^2$ 
 is defined as follows:
 \be\la{eq:Lambda}
  \Lambda_{d,e}:=\{(a,b)\,|\,
 ae+b\equiv 0\mod 
 d,\,\,\, a\ge 0,\,\, \,b\ge 0, \,\,\,
 a+b>0\}.\ee
 \elem
 \bproof
 Any derivation 
 $\alpha x\p/\p x+\beta y\p/\p y\in\ft$ 
 of bidegree $(0,0)$ is fixed by the 
 $\Ad(\TT)$-action, hence also by the 
 $G_{d,e}$-action. 
 Since $\Vect^{\mathrm c}(\AA^2)
 =\Vect^0(\AA^2)\oplus\kk\delta$, 
 where $\delta=x\p/\p x+y\p/\p y$ 
 has bidegree $(0,0)$, it suffices 
 to prove the assertion for 
 $\Vect^0(\AA^2)$.
 
 We can easily see that $\p_{-1,n}$, 
resp. $\p_{m,-1}$,
is fixed under
the ${\Ad(G_{d,e})}$-action 
if and only if 
$n=e+kd$,
resp. $m=e'+ld$, for some $k,l\ge 0$,
 if and only if 
the bidegree $(a,b)=(-1,n)$ resp. 
$(a,b)=(m,-1)$
satisfies the relation 
$ae+b\equiv 0\mod d$.
 
 Recall that the weighted 
 component 
 of $\Vect^0(\AA^2)$ 
 of bidegree $(a,b)$, 
 where  $a,b\ge 0$, 
 is one-dimensional 
 and is generated by
\[\p_{a,b}
=x^ay^b((b+1)x\p/\p x- (a+1)y\p/\p y).\] 
Since $x\p/\p x$ and $y\p/\p y$ 
are fixed under
the ${\Ad(G_{d,e})}$-action,  $\p_{a,b}$ 
is also fixed if and only if $G_{d,e}$ fixes 
$x^ay^b$, 
if and only if $ae+b\equiv 0\mod d$. 
 \eproof
 \noindent {\bf Question 2.} 
 \emph{For which pairs $(d,e)$ 
 the equality 
\be\la{eq:eq}
\Lie(\cN_{d,e})=
 \Lie({\Aut}_0(\AA^2))^{\Ad(G_{d,e})}\ee
 holds? In particular, is it true that
 \[xy(x\p/\p x-y\p/\p y)\in 
 \Lie({\Aut}_0(\AA^2))^{\Ad(G_{4,3})}
 \setminus \Lie(\cN_{4,3}),\]
 and, more generally, that 
 \eqref{eq:eq} holds if and only 
 if $ee'\in\{1, d+1\}$?} 
 
 \medskip
 
See 
 Example \ref{exa:veron}  and
 Proposition \ref{lem:5.13} below 
 for some partial results. 
\bexa
\la{exa:veron} 
Suppose that $ee'=1$, 
that is, $e=e'=1$.
In this case, 
the derivations
$\p^+_{1+kd}$ and $\p^-_{1+ld}$ 
belong to $\Lie(\Aut(X_{d,1}))$, 
see e.g. \eqref{eq-2} and \eqref{eq-3}.
Thus, $\Lie(\Aut(X_{d,1}))$ contains 
$\p^+_1=y\p/\p x$ 
and $\p^-_1=x\p/\p y$. Then also 
$[\p^-_1,\p^+_1]=
x\p/ \p x - y\p/\p y\in\Lie(\Aut(X_{d,1}))$,
and therefore 
\be\la{eq:sl2} \Lie(\Aut(X_{d,1}))\supset
\kk  y\p/\p x\oplus 
\kk x\p/ \p y\oplus \kk(x\p/ \p x - y\p/\p y)
\simeq\fs\fl(2,\kk).\ee
From \eqref{eq:comm} 
we deduce 
\[\p_{ld, kd}=[\p^-_{1+ld},\,\p^+_{1+kd}]  
\in \Lie(\Aut(X_{d,1}))
\quad \forall k,l\ge 0.\]
In particular, $\p_{ld,0}\in 
\Lie(\Aut(X_{d,1}))$ for $l\ge 0$. 
Then also 
\[\ad(\p_1^+)^s( \p_{ld,0})=
\alpha_s\p_{ld-s,s}\in \Lie(\Aut(X_{d,1})) 
\quad\forall s=1,\ldots,ld,
\,\,\,\text{where}\,\,\,\alpha_s> 0,\]
see \eqref{eq:comm}. 
Thus, $\p_{a,b}\in  \Lie(\Aut(X_{d,1}))$
for all $(a,b)$ with $a,b\ge 0$ and
$a+b\equiv 0\mod d$.
Since
$\ft_2\subset  \Lie(\Aut(X_{d,1}))$
it follows from \eqref{eq:long} that 
\[ \Lie(\Aut(X_{d,1}))=\Lie(\cN_{d,1})=
 {\Lie({\Aut}_0(\AA^2))^{\Ad(G_{d,1})}}.\]
\eexa
\brem The surface
$X_{d,1}$ is isomorphic to the $d$-th 
Veronese surface 
$V_d\subset\AA^{d+1}$, 
where $d\ge 2$, that is, to
the affine cone 
over the 
rational normal curve 
\[C_d=\{(u^d:u^{d-1}v:\ldots:v^d)\,|\,
(u:v)\in\PP^1\}\subset\PP^d.\]
Note that $V_d$ is the closure 
of an ${\SL}(2,\kk)$-orbit in
the irreducible 
${\SL}(2,\kk)$-module 
$\AA^{d+1}$
of degree $d$ binary forms, 
see e.g. \cite{Pop73}. 
The  $\fs\fl(2,\kk)$-subalgebra 
\eqref{eq:sl2}
of $\Lie(\Aut(V_d))$
is the tangent Lie algebra 
of the ${\SL}(2,\kk)$-action on $V_d$.
The ${\SL}(2,\kk)$-module 
$\AA^{d+1}$ can be considered 
as the Zariski 
tangent space of the cone $V_d$ 
at its vertex, equipped 
with the induced ${\SL}(2,\kk)$-action.
In turn, the ${\SL}(2,\kk)$-action 
on the pair $(\AA^{d+1}, V_d)$
induces the natural 
${\rm PSL}(2,\kk)$-action on the pair
$(\PP^d, C_d)$.

The interested reader can find
further information on automorphisms and 
derivations of Veronese subalgebras 
over a field of characteristic zero
 in \cite{AMLU25} and \cite{AU23}.
 \erem
\bexa\la{exa:KZ}
We have
$\Lie(\cN_{3,2})=
\Lie({\Aut}_0(\AA^2))^{\Ad(G_{3,2})}$,
see
\cite[Lemma~6.2.1]{KZ24}.
In this example, $e'=e=2$ 
and $d=ee'-1=3$.
\eexa
The last example 
can be generalized as follows.
\bprop\la{lem:5.13} The equality 
\eqref{eq:eq} holds provided that 
$ee'=d+1$\footnote{
For $d>2$ the latter equality is equivalent 
to the conditions
$e>1$ and $e | (d+1)$ (cf. 
Theorem \ref{thm:1.3}), 
and also to the conditions 
$e'>1$ and $e'|(d+1)$. }.
\eprop
\bproof 
 Recall that $\p^+_e\in \fu^+_{d,e}$ and 
 $\p^-_{e'}\in \fu^-_{d,e}$. 
 Fix an arbitrary $(a,b)\in\Lambda_{d,e}$, 
 see \eqref{eq:Lambda}. 
 We will show that
 under our assumptions,
\[\p_{a,b}\in \langle
 \p^+_e, \p^-_{e'}
\rangle_{\Lie}\subset \Lie(\Aut(X_{d,e})).\]
Since $\fu^\pm_{d,e}$ and $\ft_2$
are also subalgebras of $\Lie(\Aut(X_{d,e}))$
(see Remark \ref{rem:5.9}.2),
the above  inclusion
implies \eqref{eq:eq}. 

We have an equivalence
\[ae+b\equiv 
0\mod d\Leftrightarrow 
a+be'\equiv 0\mod d.\]
Let $u=(-1,e)$ and $v=(e',-1)$. 
There exists a pair 
$(n,m)$ of strictly positive integers 
such that
\be\la{eq:dioph} (a,b)=nu+mv
=(me'-n,ne-m).\ee
Indeed,
since by assumption $ee'=d+1$, 
the system \eqref{eq:dioph} 
of diophantine 
equations admits 
an integer solution 
\[(n,m)=\left(\frac{a+be'}{d}, 
\frac{ae+b}{d}\right),\quad
\text{where}\quad n,m>0.\]

By \eqref{eq:comm} we have 
\[u+v=(e'-1,e-1)\quad\text{and}\quad 
\p_{e'-1,e-1}=-\ad(\p^+_e)(\p^-_{e'})
\in \langle
 \p^+_e, \p^-_{e'}\rangle_{\Lie}.\]
Since $ee'=d+1$, 
the vector $(a,b)=(e'-1,e-1)$ 
is a unique lattice vector in $\Lambda_{d,e}$
satisfying the inequalities 
\be\la{eq:ineq} a<e',\quad b<e.\ee

Consider the strictly positive, 
integer valued function 
\[f\colon\Lambda_{d,e}\to d\ZZ_{>0},
\quad f(a,b)=(ae+b)+(a+be')=
(e+1)a+(e'+1)b.\] 
Suppose that \eqref{eq:ineq} 
fails for 
$(a,b)\in \Lambda_{d,e}$.
Let, for example, $a\ge e'>1$, 
and therefore
$\rho:=\lfloor \frac{a}{e'} \rfloor\ge 1$.
Consider the lattice vector
\[(a',b')=(a,b)-\rho v=
(a-\rho e',b+\rho)\in\Lambda_{d,e}.\]
Using \eqref{eq:comm} 
we deduce
\[\p_{a,b}=\beta(\ad(\p^-_{e'}))^{\rho}
(\p_{a',b'}),
\quad\text{where}\quad \beta=
[(b+2)(b+3)\cdots (b+\rho+1)]^{-1}\neq 0,\]
and 
\[ f(a',b')=f(a,b)-\rho d<f(a,b).\]
By symmetry, if $b\ge e>1$, then 
letting $\eta=\lfloor \frac{b}{e} \rfloor\ge 1$
and $ (a',b')=(a,b)-\eta u\in\Lambda_{d,e}$
 we obtain
\[\p_{a,b}=
\alpha(\ad(\p^+_e))^{\eta}(\p_{a',b'}), 
\quad\text{where}\quad 
\alpha\neq 0,\]
and
\[f(a',b')
=f(a,b)-\eta d<f(a,b).\]
Thus, any $(a,b)\in\Lambda_{d,e}$ 
distinct from 
$(e'-1,e-1)=u+v$ can be 
written in one of the forms 
\[(a,b)=(a',b')+\rho v\quad\text{or}
\quad (a,b)=(a',b')+\eta u,\]
where $\rho, \eta\ge 1$ 
and $f(a',b')<f(a,b)$.
The minimum of $f$ 
on $\Lambda_{d,e}$ 
is equal to $2d$ and is achieved at 
a single point $u+v=(e'-1,e-1)$. 
Continuing in this way, 
we finally achieve the point 
$(a',b')=(e'-1,e-1)$,
where the function $f$ 
attains its minimum. 

Summarizing,
any $(a,b)\in\Lambda_{d,e}$ 
can be written as
\[(a,b)=u+v+\rho_1v+\eta_1u+\cdots
+\eta_{k-1}u+\rho_k v,\quad k\ge 0,\]
where 
\[1+\rho_1+\cdots+\rho_k=m
\quad\text{and}\quad 
1+\eta_1+\cdots+\eta_{k-1}=n,\]
one of the $\rho_1,\rho_k$ 
or both can be zero, while the other 
coefficients are positive,
and all partial sums
\[u+v, \quad u+v+\rho_1v, 
\quad u+v+\rho_1v+
+\eta_1 u,\quad\ldots,\quad(a,b)\]
are lattice vectors from $\Lambda_{d,e}$. 
This allows to represent 
$\p_{a,b}$ as follows:
\be\la{eq:inclusion}
 \p_{a,b}=c\,(\ad(\p^-_{e'}))^{\rho_k}
(\ad(\p^+_e))^{\eta_{k-1}}\cdots
(\ad(\p^+_e))^{\eta_1}
(\ad(\p^-_{e'}))^{\rho_1}{\ad}(\p^+_e)
(\p^-_{e'})\in  \langle
 \p^+_e, \p^-_{e'}\rangle_{\Lie},\ee
 where $c\in\kk$ is nonzero, 
 see \eqref{eq:comm}.
\eproof
\bcor 
Under the assumption $ee'=d+1$ 
 we have 
 \[\Lie(\Aut(X_{d,e}))=\ft_2\oplus \langle
 \p^+_e, \p^-_{e'}\rangle_{\Lie}.\]
\ecor
\bproof According to \eqref{eq:inclusion},
$\p_{a,b}\in \langle
 \p^+_e, \p^-_{e'}\rangle_{\Lie}$ for all 
 $(a,b)\in\Lambda_{d,e}$. 
In particular, for all $k,l>0 $
we have
\be\la{eq:belongs} \p_{0,kd}, \,\,\p_{ld,0}\in\langle
 \p^+_e, \p^-_{e'}\rangle_{\Lie}. \ee
 Moreover, by  \eqref{eq:comm},
\[\p^+_{e+kd}=(e+kd+1)^{-1}
\ad(\p^+_e)(\p_{0,kd})
\quad\text{and}\quad
 \p^-_{e'+ld}=-(e'+ld+1)^{-1}
\ad(\p^-_{e'})(\p_{ld,0})\] 
also belong to
$ \langle\p^+_e, \p^-_{e'}\rangle_{\Lie}$.
 Now the claim follows. 
\eproof
 By combining Theorem \ref{thm:1.2}, 
 Theorem \ref{main}, Lemma \ref{lem:5.5},
 and  Corollary \ref{cor:5.8}, we arrive 
 at the following conclusion.
 \bcor\la{cor:new}  The locally
integrable Borel subalgebras of 
$\Lie(\Aut(X_{d,e}))$ are precisely 
the $\Ad$-conjugates 
of $\ft_2 \ltimes \fu^{\pm}_{d,e}$.
 \ecor
 \brem\la{rem:ref-question}
 However,  $\ft_2 \ltimes \fu^{\pm}_{d,e}$ 
 is not always a Borel subalgebra of 
 $\Lie(\Aut(X_{d,e}))\simeq \Lie(\cN_{d,e})$.
 Indeed, suppose that $ee'=d+1$. 
 Then we have $\p_{0,d}\in\Lie(\Aut(X_{d,e}))$,
 see \eqref{eq:belongs}. In this case,
the solvable Lie subalgebra 
\[\fh:= \langle \ft_2 \ltimes \fu^{\pm}_{d,e}, 
\p_{0,d}\rangle_{\Lie}\subset\Lie(\Aut(X_{d,e}))\]
strictly contains $\ft_2 \ltimes \fu^{\pm}_{d,e}$;
see the proof of Proposition \ref{prop:non-borel}.
 \erem
\subsection{The Lie algebra 
of zero divergence vector fields 
on $X_{d,e}$ is not simple} 
\la{ss:6.3}
Recall that the subalgebra 
$\Vect^0(\AA^2)$
of divergence-free vector fields 
is the unique proper ideal of 
the Lie algebra
$\Lie(\Aut(\AA^2))=\Vect^{\mathrm c}(\AA^2)$.
Moreover, the codimension one 
subalgebra ${\Vect^0(\AA^2)\subset
\Lie(\Aut(\AA^2))}$
is simple, 
see \cite[Lemma~ 3]{Sha81}.

In the Lie algebra $\Lie(\cN_{d,e})$,
we consider the codimension 
one ideal
\[\fg_{d,e}=\Lie(\cN_{d,e}) 
\cap \Lie(\mathfrak{U}_2)=
\Lie(\cN_{d,e})
 \cap{\Vect}^0(\AA^2).\] 
Since $\Lie(\Aut(X_{d,e}))$ 
can be naturally identified 
with $\Lie(\cN_{d,e})$, see
Corollary \ref{cor:5.8},
we can consider $\fg_{d,e}$ 
as a maximal ideal of 
$\Lie(\Aut(X_{d,e}))$.
\bprop\la{prop:nonsimple} 
\la{prop:non-simple}
The Lie  subalgebra 
$\fg_{d,e}\subset \Lie(\Aut(X_{d,e}))$
 is not simple.
\eprop
\bproof
 In the lattice 
 \[L_{d,e}=\{(a,b)\in\ZZ^2\, |\, 
 ae+b\equiv 0 \mod d\}\]
we consider the submonoid
\[\hat\Lambda_{d,e} = 
\{(a,b)\in L_{d,e}\,|\,a,b \ge -1,
 \quad (a,b) \neq (-1,-1)\}.\]
 It properly contains the 
 submonoid $\Lambda_{d,e}$ 
 of \eqref{eq:Lambda}. More precisely,
 \[\hat\Lambda_{d,e}=\Lambda_{d,e}
 \cup \{(-1,e+kd)\,|\,k\ge 0\}
 \cup\{(e'+ld,-1)\,|\,l\ge 0\}\cup\{(0,0)\}.\]
 We also consider the submonoid 
 \[\lambda_{d,e}:=\hat\Lambda_{d,e}
 \setminus\{(-1,e),(e',-1), (0,0)\}\]
and the homogeneous 
subalgebra $I_{d,e}\subset 
\Lie({\Aut}_0(\AA^2))^{\Ad(G_{d,e})}$ 
of codimension 4
supported by $\lambda_{d,e}$, where
\[I_{d,e}=\bigoplus_{k\ge 1} 
\kk\p^+_{e+kd}
\oplus\bigoplus_{l\ge 1} \kk\p^-_{e'+ld}
\oplus\bigoplus_{(a,b)\in
\Lambda_{d,e}} \kk\p_{a,b},\]
cf. \eqref{eq:long}. 
We claim that 
$I_{d,e}$ is an ideal of 
$\Lie({\Aut}_0(\AA^2))^{\Ad(G_{d,e})}$. 
Using 
\eqref{eq:comm},
it suffices to check the inclusion
\[(\hat\Lambda_{d,e} + \lambda_{d,e}) 
\cap \hat\Lambda_{d,e} \subset \lambda_{d,e}.\]
Suppose this inclusion 
is not true. Then 
the left-hand side contains 
one of the lattice vectors
$(-1,e), \,(e',-1)$ and $(0,0)$. 
Assuming that $(-1,e)=\hat v+v$, 
where $\hat v\in\hat\Lambda_{d,e}$ 
and $v\in\lambda_{d,e}$, 
we obtain 
\[(-1,e)=(-1, e+kd)+(0,ld)\quad
\text{for some}\quad k,l>0,\] 
which is impossible.
The same argument applies to $(e',-1)$.
Similarly, if $(0,0)=\hat v+v$, 
then either $\hat v=v=(0,0)$, or
$\{\hat v,v\}=\{(-1,1),\,(1,-1)\}$. But
this is impossible because
$(0,0),\,(-1,1),\,(1,-1)\notin\lambda_{d,e}$.

Since $I_{d,e}$
 is a 
 homogeneous ideal 
of codimension $4$
 of $\Lie({\Aut}_0
 (\AA^2))^{\Ad(G_{d,e})}$,
 then also 
$I^0_{d,e}:=I_{d,e}\cap\fg_{d,e}$ 
is a homogeneous ideal 
of codimension $3$
of the subalgebra
$\fg_{d,e}\subset \Lie({\Aut}_0
 (\AA^2))^{\Ad(G_{d,e})}$.
Indeed, it is transversal to 
the vector subspace
\[\kk\p^+_{e}\oplus\kk\p^-_{e'}
\oplus\kk\p_{0,0}\subset\fg_{d,e}.\]
Therefore,
the Lie algebra 
$\fg_{d,e}$ is not simple. 
\eproof
\brems $\,$

1.
 This proof does 
not work for 
$\Lie(\mathfrak{U}_2)$, 
because this Lie algebra contains 
the derivations $\p/\p x$ and $\p/\p y$
whose adjoint action corresponds to 
the left or down shift
by one unit
on the lattice $\ZZ^2$, 
see \eqref{eq:comm}. 
And, indeed, this algebra is simple, 
see \cite[Lemma~3]{Sha81}.

2. Any normal affine 
toric surface $X$ carrying no 
non-constant invertible regular 
function is isomorphic to either $\AA^2$ 
or to some $X_{d,e}=\AA^2/\G_{d,e}$. 
It follows from 
Proposition \ref{prop:nonsimple}
that $X\simeq\AA^2$ 
(that is, $X$ is smooth) 
if and only if the maximal ideal 
of divergence zero 
vector fields from $\Lie(\Aut(X))$ 
is simple 
(cf. \cite{Sie92}). 
\erems
{\bf Acknowledgments}. 
It is a pleasure to thank 
Hanspeter Kraft
for his assistance and for 
informing us in
important results from 
the preprint
\cite{CKRVS26} used in the proofs.
We would like to acknowledge
 the remarkable work 
of the anonymous reviewer,
whose careful reading, comments,
and suggestions helped to
avoid several inaccuracies 
and enriched the paper. In particular,
Proposition \ref{prop:non-borel}, 
Corollary \ref{cor:new}, as well as
Remarks \ref{rem:suggest} and 
 \ref{rem:ref-question},
were added following the reviewer 
suggestions. 

\end{document}